# Airline recovery problem under disruptions: A review

Shuai Wu[1], Enze Liu[1], Rui Cao[1], Qiang Bai[1*]

[1] School of Transportation Engineering, Chang'an University, Xi'an, Shaanxi, China, 710064

* Corresponding author.

**Abstract:** In practice, both passenger and cargo flights are vulnerable to unexpected factors, such as adverse weather, airport flow control, crew absence, unexpected aircraft maintenance, and pandemic, which can cause disruptions in flight schedules. Thus, managers need to reallocate relevant resources to ensure that the airport can return to normal operations on the basis of minimum cost, which is the airline recovery problem. Airline recovery is an active research area, with a lot of publications in recent years. To better summarize the progress of airline recovery, first of all, keywords are chosen to search the relevant studies, then software is used to analyze the existing studies in terms of the number of papers, keywords, and sources. Secondly, the airline recovery problem is divided into two categories, namely Passenger-Oriented Airline Recovery Problem (POARP) and Cargo-Oriented Airline Recovery Problem (COARP). In POARP, the existing studies are classified according to recovery strategies, including common recovery strategies, cruise speed control strategy, flexible aircraft maintenance strategy, multi-modal transportation strategy, passenger-centric recovery strategy, and clubbing of flights strategy. Moreover, the POARP is discussed from the perspectives of disruptions, problem types, objective functions, and solution methods. Thirdly, POARP and COARP are compared from the perspectives of timeliness, subjectivity, flexibility, transferability, and combinability. Finally, the conclusions are drawn and future study directions are provided. For future studies, it is recommended to conduct more in-depth research on dynamic and real-time recovery, incorporating human factors into the modeling, multi-modal transportation coupling, optimization of other airport processes, combination of robust scheduling and airline recovery, and optimization algorithm improvement.

**Key words:** airline recovery problem; passenger-oriented recovery; cargo-oriented recovery; disruption; recovery strategy

# 1 Introduction

Air transportation is one of major transportation mode and plays an important role in the economy development. Before the COVID-19 pandemic, there were more than 4.5 billion air passengers and more than $6.7 trillion cargo carried by air transportation globally in 2019 (de Juniac, 2020). Proactive strategies, such as forecasting air passenger and cargo demand, planning flight schedules, fleet assignment, and aircraft routing, are typically developed to ensure the normal operation of the air transportation system (Derigs and Friederichs, 2012; Zhou et al., 2020).

However, flights are easily disrupted by unexpected events, such as adverse weather, airport flow control, crew absence, unexpected aircraft maintenance, and epidemic (Hsu and Shih, 2010; Evans and Schäfer, 2014; Chen and Chou, 2017; Glass et al., 2022; Wen et al., 2022). Proactive strategies cannot guarantee a good response to disruptions in these events and may lead to flight delays and even large-scale flight cancellations (Clausen et al., 2010). Hence, reactive strategies, known as airline recovery, are needed to deal with the disruptions and restore normal operations. Reactive strategies are important for local governments, airports, airline companies, air passengers, and air cargo customers. For local governments, efficient airline recovery can enhance airport resilience, attract more passengers, and stimulate local economic growth. For airports, airline recovery can effectively reduce congestion and ensure that the airport can return to normal operation in a short time. For airlines, airline recovery can reduce costs and ensure that passengers arrive at their destinations on time, and ensure timely delivery of cargo. For air passengers, airline recovery can enable them to reach their destination as early as possible, reducing anxiety and fostering loyalty to airlines. For air cargo customers, airline recovery can ensure the timely delivery of their cargo, enhancing their satisfaction with airlines.

There is a significant number of studies on airline recovery problems. Several review papers have summarized the studies on passenger-oriented airline recovery problems (POARP) (Kohl et al., 2007; Clausen et al., 2010; Hassan et al., 2021; Su et al., 2021; Santana et al., 2023). However, these papers introduced airline recovery studies either based on the categories of models and solution methods or based on the categories of individual and integrated recovery mode. Moreover, there is a need to include the studies on cargo-oriented airline recovery problems (COARP) and systematically summarize the conclusions. Therefore, this paper categorizes the airline recovery problem into two main types: POARP and COARP. In POARP, recovery strategies are classified according to the types of strategies, including common recovery strategies, cruise speed control strategy, flexible aircraft maintenance strategy, multi-modal transportation strategy, passenger-centric recovery strategy, and clubbing of flights strategy. Furthermore, the distinctions between POARP and COARP in terms of timeliness, subjectivity, flexibility, transferability, and combinability are discussed. The purpose is to comprehensively and systematically summarize existing studies in the field of airline recovery from multiple perspectives and provide directions for future studies.

The rest of the paper is organized as follows. In Section 2, a systematic analysis of the searched literature is conducted, including the number of papers, keywords, and sources. Section 3 provides a comprehensive discussion of studies on airline recovery problems, with a primary focus on two major categories: POARP and COARP, and offers corresponding comparative analyses for these two types of studies. Finally, Section 4 presents the conclusions of this paper and outlines directions for future studies.

# 2 Literature analysis

A literature search was conducted in the Web of Science Core Collection Database using keywords "airline recovery", "air cargo recovery", and "flight recovery" with a search period up to December 2023. After conducting preliminary screening and excluding irrelevant papers, a total of 1054 papers that align with the topic of this study were obtained. To effectively summarize the research progress on airline recovery problems, the powerful visualization capabilities of VOSviewer software (van Eck and Waltman, 2010) was used to visually organize and analyze the literature from three aspects: number of papers, keywords, and sources.

## 2.1 Number of papers

In the research area of airline recovery, Fig.1 presents the number of papers published each year from 1984 to 2023. It can be seen that, starting from 2010, the number of papers has been increasing sharply, indicating that the problem is becoming a popular reach field. As a result, there is a growing need to summarize and synthesize existing studies and to identify future study directions.

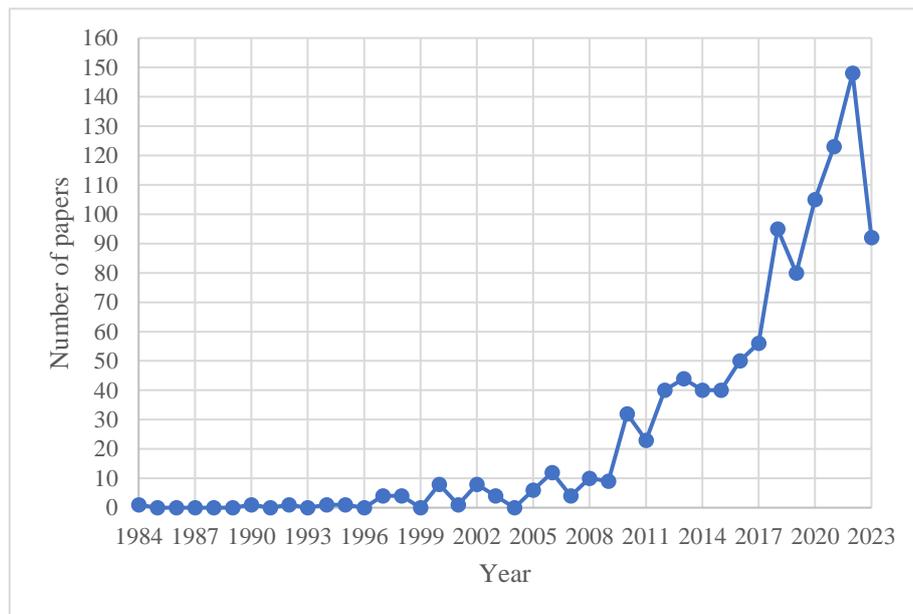

**Fig.1 The number of published papers from 1984 to 2023**

The VOSviewer is used to visualize the number of publications and citations from different countries and organizations. The result is shown in Figs.2-4. It is seen that China has the most publications, followed by the countries with renowned transportation universities, such as the United States and England. In terms of publications by organizations, traditional transportation-focused universities rank prominently. Chinese universities, including Shanghai Jiao Tong University, Beijing Jiaotong University, and Hong Kong Polytechnic University, have a substantial publication output. American institutions like MIT and University of Michigan are indeed prominent in the field of airline recovery, along with other internationally renowned universities such as Delft University of Technology and Imperial College London, which have made significant contributions to this research field. Lastly, regarding the countries that cite these publications, China, the United States, and England remain at the forefront, which indicates that scholars from these countries are actively advancing research on airline recovery problem. In summary, scholars from traditional transportation-focused organizations naturally show a higher interest in this field. China, the United States, and England are among the top countries in terms of publication output. Scholars from other countries have also made significant contributions to research in this area. Overall, it is seen that airline recovery problem has become a hot topic of study for scholars worldwide.

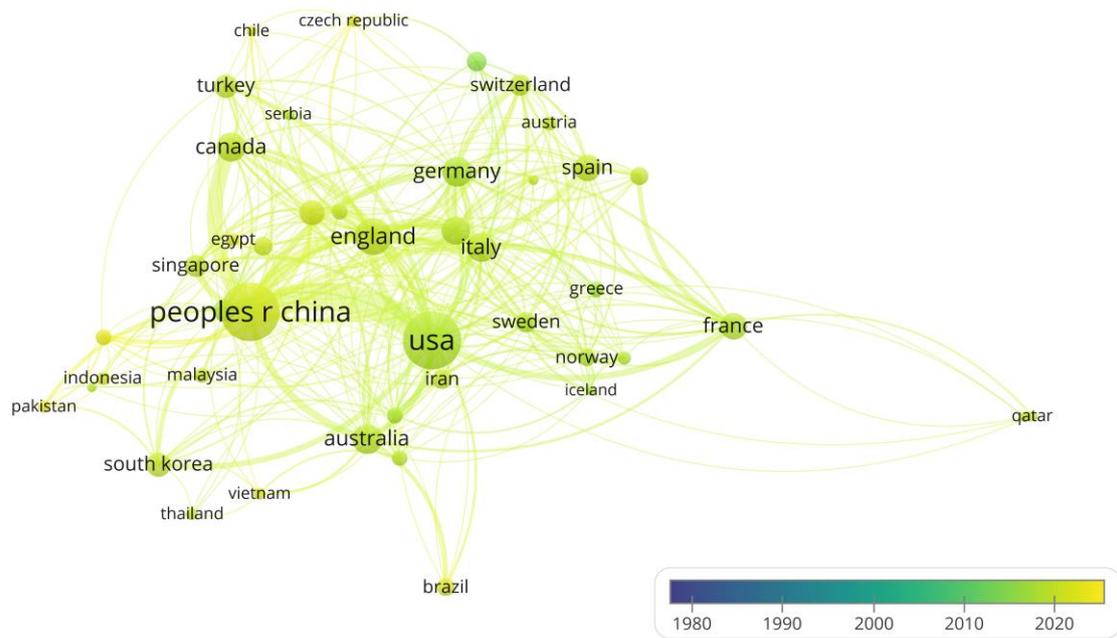

**Fig.2 Mapping knowledge domain of published papers by countries**

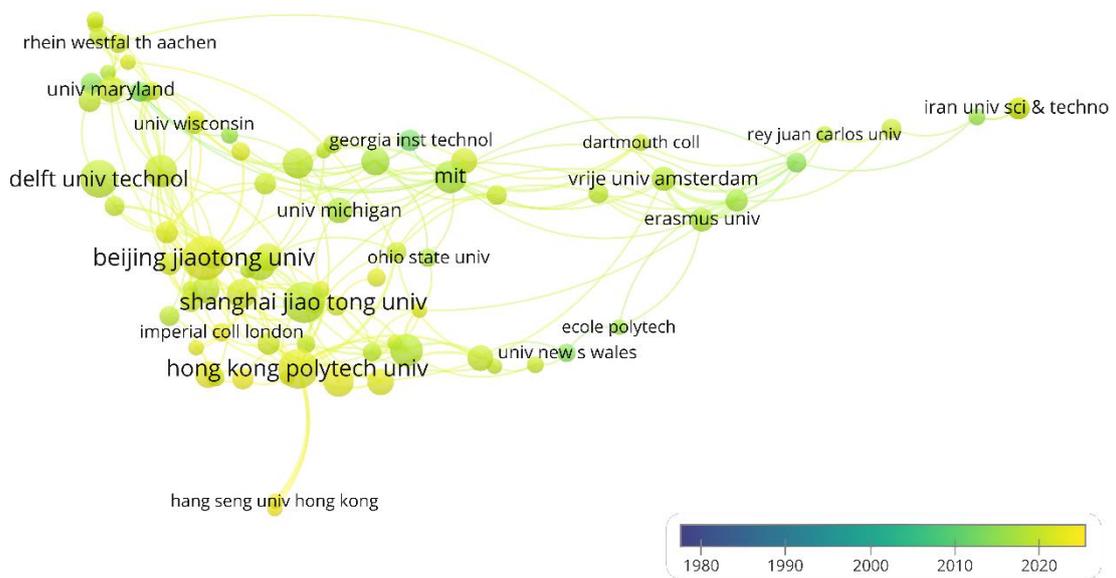

**Fig.3 Mapping knowledge domain of published papers by organizations**

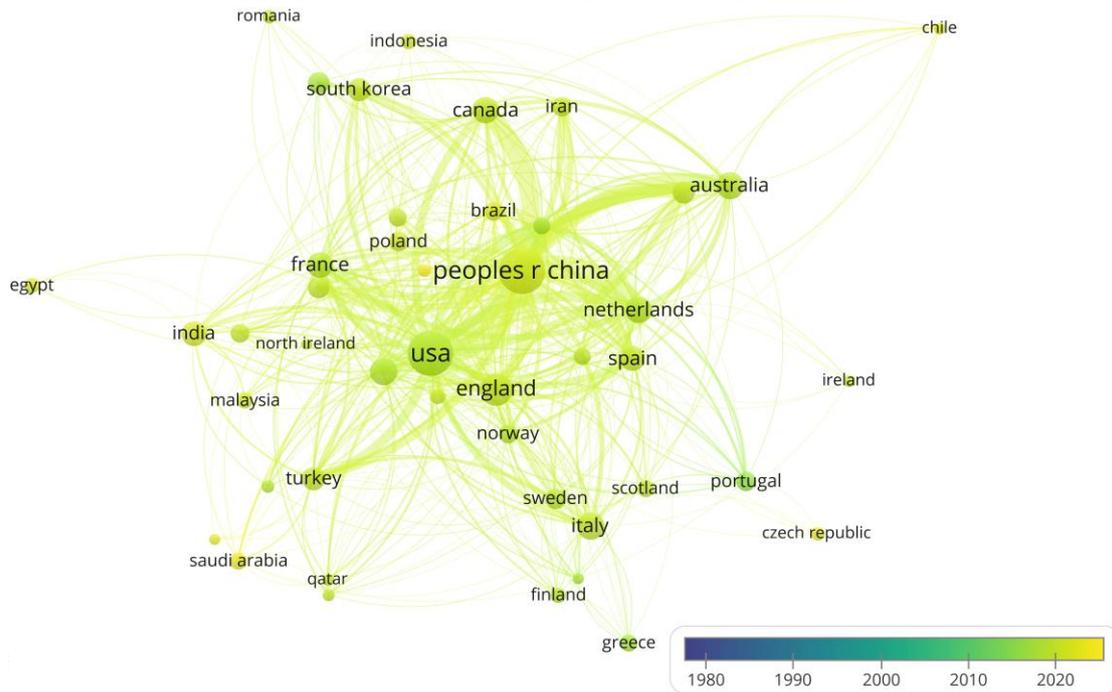

**Fig.4 Mapping knowledge domain of total citations by countries**

## 2.2 Keywords

To gain a better understanding of key subareas in this field, VOSviewer is used to visually present the keywords in the existing studies, as shown in Fig.5. It can be seen that that "recovery", "disruption management", and "irregular operations" are primary research subareas. Keywords such as "airport closure", "hub closures", "flight cancellations", and "delay propagation" highlight the most commonly studied disruption cases. Keywords such as "aircraft", "flight", "crew", "passenger", and "air cargo" show the main subjects of studies in airline recovery problem. Keywords such as "integer programming", "optimization model", and "optimization" indicate that the airline recovery is an optimization problem. In terms of solution algorithms, "heuristic algorithm" and "genetic algorithm" are popular methods in the literature. The keywords co-occurrence mapping knowledge domain effectively shows the primary study contents and methods in this field, enabling readers to quickly understand the core aspects of this research field.

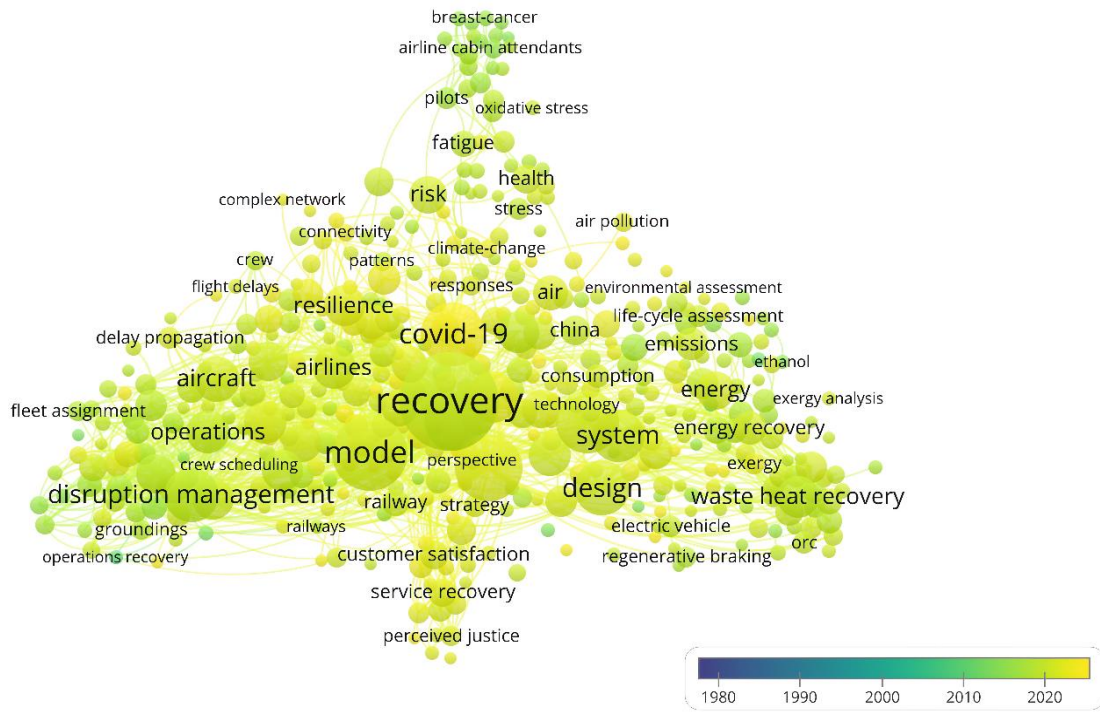

**Fig.5 Keywords co-occurrence mapping knowledge domain**

## 2.3 Sources

The numbers of published papers on airline recovery problem in various journals are shown in Table 1. Firstly, airline recovery falls within the domain of air transportation problem, requiring the management of disrupted flights. Therefore, *Journal of Air Transport Management* has the highest number of papers in this area. Secondly, traditional transportation journals such as *Transportation Research Record*, *Transportation Research Part C*, *Transportation Science*, *IEEE Transactions on Intelligent Transportation Systems*, *Transport Policy*, *Transportation Research Part E*, *Transportation Research Part A*, and *Transportation Research Part B* all have relevant research papers. Furthermore, airline recovery problem involves operations research and algorithm-related contents, resulting in a substantial number of papers in journals within the computer science and operations research domain, such as *Computers & Operations Research*, *Computers & Industrial Engineering*, and *European Journal of Operational Research*. In summary, airline recovery problem has become an interdisciplinary problem that involves transportation, management, computer science, and operations research, making it a subject worthy of continued exploration by scholars.

**Table 1 Journal ranking by number of publications on airline recovery from 1984 to 2023**

| Rank | Journals | Number of Publications |
| --- | --- | --- |
| 1 | Journal of Air Transport Management | 37 |
| 2 | Transportation Research Record | 35 |
|   | Sustainability |   |
| 3 | Transportation Research Part C: Emerging Technologies | 26 |
| 4 | Transportation Science | 20 |

| 5 | IEEE Transactions on Intelligent Transportation Systems | 18 |
| | Transport Policy | |
| 6 | Computers & Operations Research | 17 |
| 7 | Transportation Research Part E: Logistics and Transportation Review | 16 |
| | Computers & Industrial Engineering | |
| 8 | Journal of Advanced Transportation | 15 |
| 9 | European Journal of Operational Research | 14 |
| | Transportation Research Part A: Policy and Practice | |
| 10 | Transportation Research Part B: Methodological | 13 |

# 3 Airline recovery problem

The review is conducted by categorizing existing studies on airline recovery into POARP and COARP. Within the POARP category, the literature is further classified according to the types of recovery strategies employed. Section 3.1 summarizes existing studies on POARP; Section 3.2 provides an overview of COARP; and Section 3.3 summarizes the similarities and differences between POARP and COARP.

## 3.1 Passenger-oriented airline recovery

Air passenger-oriented airline recovery strategies can be grouped into common strategies which are commonly used in most cases and special strategies which are typically used in special cases. Details of these strategies are presented in Table 2.

**Table 2 Summary of recovery strategies on POARP**

| Categories | Contents |
|---|---|
| Common strategies | Flight delay, flight cancellation, aircraft swap, ferry aircraft, crew delay, crew cancellation, crew swap, crew deadhead, use reserve aircraft and crews |
| Special strategies | Aircraft cruise speed control strategy |
| | Flexible aircraft maintenance strategy |
| | Multi-modal transportation strategy |
| | Passenger-centric recovery strategy |
| | Clubbing of flights strategy |

### 3.1.1 Common recovery strategies

In practice, various disruptions occur, including adverse weather, flight delays, airport flow control, and unexpected aircraft maintenance. To facilitate studies, most literature focuses on specific disruption cases such as flight delays and flight cancellations. For these disruptions, several commonly-used recovery strategies can be adopted to minimize the impacts on passengers and airline companies, including flight delay and flight cancellation, aircraft swap, using reserve aircraft and crews, ferrying aircraft, and crew deadhead. Flight delay and flight cancellation strategies are straightforward and do not require further explanation. Using reserve aircraft or crew means that airports allocate additional resources to deal with the disruptions. The ferrying aircraft strategy involves transferring an idle aircraft from another location to the designated airport to operate a flight. However, this strategy comes with higher costs, so airlines use it sparingly. Crew deadhead is similar to ferrying aircraft, where crews travel on an aircraft to reach a specified airport to perform their assigned flight duties. Furthermore, the basic principle of the aircraft swap strategy is illustrated in Fig.6.

There are two same types of aircraft operating flight tasks. Aircraft $a_1$ is paired with crew $c$ to operate flights $f_1$ and $f_2$; and aircraft $a_2$ is paired with crew $d$ to operate flights $f_3$ and $f_4$. Suppose an unexpected event causes $f_3$ to be delayed in arriving at the airport. If the aircraft swap strategy is not adopted, this delay will propagate to $f_4$. However, by executing the aircraft swap strategy, after $f_1$ arrives at the airport, $a_1$ with crew $c$ can operate $f_4$; and after $f_3$ arrives at the airport, $a_2$ with crew $d$ can operate $f_2$. This can reduce delay caused by unexpected events, improve airport efficiency, and ensure passenger satisfaction.

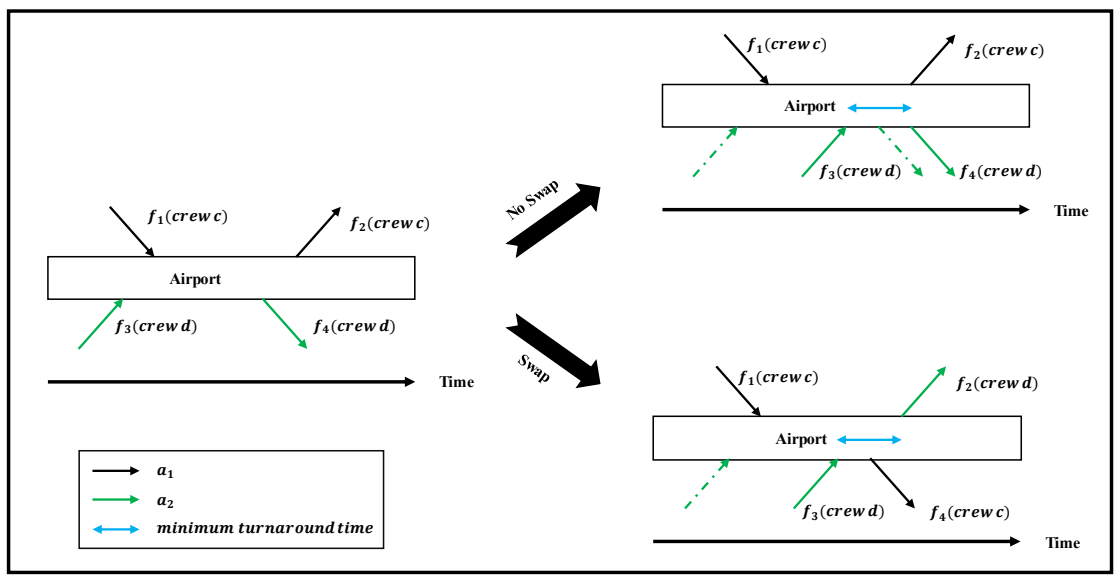

**Fig.6 Aircraft swap strategy**

Based on the types of airline recovery strategies, the rest of this subsection is divided into six parts to provide a clearer summary of existing academic achievements. Each part discusses one of the commonly-used airline recovery strategies, including aircraft recovery, crew recovery, passenger recovery, aircraft and passenger recovery, aircraft and crew recovery, and integrated airline recovery (Hassan et al., 2021).

**3.1.1.1 Aircraft recovery**

The Aircraft Recovery Problem (ARP) has garnered early attention from scholars due to the criticality and high cost associated with aircraft, making it one of the most vital resources at an airport. A summary of existing studies can be found in Table 3.

One of the most studied two disruptions is Aircraft-on-Ground (AOG), as can be seen from Table 3. Teodorović and Guberinić (1984) first studied ARP for the AOG. They developed a model with the objective of minimizing the total delay time using the flight delay strategy and designed a Branch-and-Bound (B&B) algorithm to solve this problem. Teodorović and Stojković (1990) further incorporated the flight cancellation strategy into the model, aiming to minimize both delay time and the number of cancelled flights. Their study framework provided valuable insights for subsequent studies, which mainly focused on recovery strategies, objective functions, and solution methods. Firstly, in terms of recovery strategies, aircraft swap was developed to restore the airport operation (Hu et al., 2017; Zhao et al., 2023). Then, in terms of objective functions, several studies set minimizing Total Recovery Cost (TRC) as the objective function, which includes delay cost, cancellation cost, and swap cost (Zhu et al., 2015; Zhao et al., 2023). Lastly, many effective and efficient solution methods were developed to solve the ARP, such as Greedy Simulated Annealing (GSA) algorithm (Zhu et al., 2015), Binary Search-Minimum Cost Flow (BSMCF) algorithm and Combined ε-constraints method and Neighborhood Search (CNS) (Hu et al., 2017), and Rolling Horizon (RH) algorithm (Zhao et al., 2023).

The other most studied two disruption is Airport Closure (AC). Thengvall et al. (2000) are indeed the pioneers in studying ARP in AC. They developed a model with an objective of maximizing profit and used CPLEX to solve it. In a later study, Thengvall et al. (2001) expanded their work by considering more recovery strategies and formulating ARP as a multi-commodity network flow model. Thengvall et al. (2003) designed the Bundle Algorithm (BA), which significantly accelerated the running time compared to CPLEX. Different from the above studies, subsequent studies mainly focused on the objective functions of minimizing TRC (Le et al., 2013) or delay time (DT) (Liu et al., 2008, 2010; Lee et al., 2022). Le et al. (2013) designed a Genetic Algorithm (GA) to minimize TRC. To minimize DT, Liu et al. (2008, 2010) developed the Method of inequalities based Multiobjective Genetic Algorithm (MMGA) to enhance the solution quality. To further improve the solution quality and reduce running time, Lee et al. (2022) innovatively adopted Reinforcement Learning (RL) to solve the model. They compared the proposed method with solutions obtained from the Time-Line Network (TLN) model, the MMGA proposed by Liu et al. (2008), and WSAP (Without Swapping) scheme. Results indicated that RL could acquire high-quality solutions in a relatively short time, validating the effectiveness of the method.

Several studies have explored both AOG and AC disruptions (Zhao and Guo, 2012a, 2012b; Arias et al., 2013; Zhang, 2017; Wang et al., 2019; Ji et al., 2021; Liu ,2022). The common objective function across these studies is to minimize the total recovery cost. The innovation lies mainly in the improvement of solution methods, such as Greedy Random Adaptive Search Procedure-Ant Colony Optimization (GRASP-ACO) algorithm (Zhao and Guo, 2012b) and Two-Stage Heuristic (TSH) algorithm (Zhang, 2017). Some studies directly considered the quantity of delayed or cancelled flights as disruption cases for studying. For instance, Zhao and Guo (2012a) developed an improved GRASP algorithm; Arias et al. (2013) innovatively combined simulation and optimization algorithms to study ARP; Ji et al. (2021), referring to the work of Gao et al. (2012), integrated considerations of flight priority and airport capacity limitations, and designed TSH algorithm to address these combined factors. Wang et al. (2019) used simulation to explore recovery schemes under various adverse weather conditions. Liu (2022) conducted a study on ARP in the context of flight delays caused by a snowstorm.

From Table 3, it can also be seen that flight delay and flight cancellation are the most adopted strategies and that minimizing the DT is the most used objective function. Regarding the solution methods, heuristic algorithms dominate other methods. Overall, common recovery strategies for the ARP have been well studied.

Table 3 Summary of aircraft recovery studies using common recovery strategies

| Studies | Disruptions | | | | | Recovery strategies | | | | | | Objective functions Min | | | Max Profit | Solution methods/ tools |
|---|---|---|---|---|---|---|---|---|---|---|---|---|---|---|---|---|
| | AOG | AC | FD | FC | AW | FD | FC | SW | Ferry | DT | TRC | NDF | NAS | Others | | |
| Teodorović and Gubernić (1984) | Y | | | | | Y | | | | Y | | | | | | B&B |
| Teodorović and Stojković (1990) | Y | | | | | Y | Y | | | Y | | Y | | | | DP |
| Thengvall et al. (2000) | | Y | | | | Y | Y | Y | | | | | | | Y | CPLEX |
| Thengvall et al. (2001) | | Y | | | | Y | Y | | | | | | | | Y | CPLEX |
| Thengvall et al. (2003) | | Y | | | | Y | Y | | | | | | | | Y | BA |
| Liu et al. (2008) | | Y | | | | Y | Y | Y | | Y | | Y | Y | Y | | MMGA |
| Liu et al. (2010) | | Y | | | | Y | Y | Y | | Y | | Y | Y | Y | | MMGA |
| Zhao and Guo (2012a) | | | Y | Y | | Y | Y | | | | Y | | | | | Improved GRASP |
| Zhao and Guo (2012b) | Y | Y | | | | Y | Y | | | | Y | | | | | GRASP-ACO |
| Arias et al. (2013) | | | Y | | | Y | Y | Y | | Y | | Y | | | | Simulation |
| Le et al. (2013) | | Y | | | | Y | | | | | Y | | | | | GA |
| Zhu et al. (2015) | Y | | | | | Y | Y | | | | Y | | | | | GSA |
| Hu et al. (2017) | Y | | | | | Y | Y | Y | | Y | | | Y | Y | | BSMCF and CNS |
| Zhang (2017) | Y | Y | | | | Y | Y | | | | Y | | | | | TSH |
| Wang et al. (2019) | | | | | Y | Y | Y | | | | | | | | | Simulation |
| Ji et al. (2021) | | | Y | | | Y | | Y | | Y | | | | | | TSH |
| Lee et al. (2022) | | Y | | | | Y | | Y | | Y | | Y | | | | QL and DQL |
| Liu (2022) | | | | Y | | Y | | | | | Y | | | | | B&B |
| Zhao et al. (2023) | Y | | | | | Y | Y | Y | | | Y | | | | | RH |

Notes: AOG: aircraft-on-ground; AC: airport closure; FD: flight delay; FC: flight cancellation; AW: adverse weather; SW: aircraft swap; DT: delay time; TRC: total recovery cost; NDF: the number of disrupted flights; NAS: the number of aircraft swaps; B&B: branch-and-bound; DP: dynamic programming; BA: bundle algorithm; MMGA: the method of inequalities based multiobjective genetic algorithm; GRASP-ACO: greedy random adaptive search procedure - ant colony optimization; GA: genetic algorithm; GSA: greedy simulated annealing; BSMCF: binary search-minimum cost flow; CNS: combined ε-constraints method and neighborhood search; TSH: two-stage heuristic; QL: Q-learning; DQL: double Q-learning; RH: rolling horizon; Y: included or mentioned.

**3.1.1.2 Crew recovery**

Effective scheduling of flight crew members, who are integral to carrying out flight operations, is crucial for recovery from unforeseeable disruptions. Therefore, disruption management for crew has drawn the attention of scholars all over the world (Wen et al., 2021). Wei and Yu (1997) were the first to study Crew Recovery Problem (CRP) and innovatively provided airlines with a comprehensive framework for crew management under the disruptions of flight delays and flight cancellations. Based on the framework, several studies enriched the theory (Guo et al., 2005; Chang, 2012; Liu et al., 2013). In addition, some studies explored CRP under the disruption of AOG and adverse weather (Lettovský et al., 2000; Le and Sun, 2011). Some other studies focused on the disruption of crew absence. Abdelghany et al. (2004) developed a proactive crew recovery decision support tool with the objective of minimizing recovery costs, using the RH algorithm to solve the model. Chen and Chou (2017) proposed a model with the objectives of minimizing the number of crews whose duties are changed, the number of duties assigned for different tasks, the maximal number of duty changes for a crew, the largest changed flight time for a crew, the deviation of changed duties and flight time of crews, and the number of replaced nice-to-have tasks. They combined the Nondominated Sorting Genetic Algorithm II (NSGA-II) with a Constraint-Loosening Mechanism (CLM) to solve this complex model. Schrotenboer et al. (2023) examined the impact of using reserve crews on the robustness of schedule results, modelled the planning of reserve crews using a Markov chain, and used a Branch-and-Price (B&P) algorithm to solve the model. Their proposed method aimed to create a more robust schedule, resulting in fewer crew swaps and lower recovery costs.

A summary of existing studies CRP is shown in Table 4. It can be seen that crew absence and flight delay are the most studied disruptions in CRP. In terms of recovery strategies, crew delay, deadhead, and use of reserve crews are the most commonly-used strategies, followed by crew swap. Regarding objective functions, the total recovery cost is the most adopted; and a few studies used the number of crew tasks and disrupted crews (Chang, 2012; Chen and Chou, 2017).

Table 4 Summary of crew recovery studies using common recovery strategies

| Studies | Disruptions | | | | | | Recovery strategies | | | | | | Objective functions Min | | | Solution methods/ tools |
|---|---|---|---|---|---|---|---|---|---|---|---|---|---|---|---|---|
| | AOG | FD | FC | AW | CA | Others | CD | CC | CS | DH | URC | TRC | NDC | NDT | Others | |
| Wei and Yu (1997) | | Y | Y | | | | Y | Y | Y | Y | Y | Y | | | | B&B |
| Lettovský et al. (2000) | Y | | | Y | | | Y | Y | Y | Y | | Y | | | | PSSM |
| Abdelghany et al. (2004) | | | | | Y | | Y | | Y | Y | Y | Y | | | | A rolling approach |
| Guo et al. (2005) | | Y | | | | | Y | | Y | Y | Y | Y | | | | Improved GA |
| Le and Sun (2011) | Y | | | Y | Y | Y | Y | Y | Y | | | Y | | | | LINGO |
| Chang (2012) | | Y | | | | | Y | | | Y | Y | | Y | Y | Y | Improved GA |
| Liu et al. (2013) | | Y | Y | | | | Y | | | Y | | | | | Y | SA |
| Chen and Chou (2017) | | | | | Y | | Y | | Y | Y | Y | | Y | Y | Y | Improved NSGA-II |
| Schrotenboer et al. (2023) | | | | | Y | | | Y | | Y | Y | Y | | | | B&P |

Notes: CA: crew absence; CD: crew delay; CC: crew cancellation; CS: crew swap; DH: deadhead; URC: use reserve crews; NDC: the number of disrupted crews; NDT: the number of disrupted tasks; PSSM: primal-dual subproblem simplex method; SA: simulated annealing; NSGA-II: nondominated sorting genetic algorithm II; B&P: branch-and-price.

### 3.1.1.3 Passenger recovery

Satisfaction is indeed the most direct evaluation criteria that passengers have of travel services. When passengers are satisfied with a particular mode of transportation, their willingness to use that mode again in the future and their likelihood to recommend it to others also tend to increase correspondingly (Yuan et al., 2021a, 2021b). For passengers, the recovery strategy implemented in disruptions is crucial. Effectively appeasing passengers' emotions is an important measure in reducing long-term losses. This is known as Passenger Recovery Problem (PRP). However, PRP is seldom studied alone. Only three relevant papers were found on this topic, as listed in Table 5. Acuna-Agost et al. (2015) studied the PRP in disruption cases involving AC, AOG, and flight delays. They developed a model with the objective of minimizing the total cost for both split and non-split groups of passengers, and designed a Network Pruning Approach (NPA) to solve it. McCarty and Cohn (2018) developed a two-stage stochastic model to deal with the rerouting of passengers. The first stage involved preemptive decisions, allocating passengers to alternative itineraries to assess the impact of the delay; the second stage aimed to further modify the itineraries to identify passengers who would still miss their flights after the delay, with the objective of minimizing delay time and recovery costs. A Benders Decomposition (BD) algorithm was used to solve the model. Wang et al. (2020) considered multi-flights and seat-sharing constraints, and formulated a model with the objective of minimizing both passenger satisfaction costs and airline costs. To solve this problem, they designed a B&P based on Column Generation (CG) framework. Results indicated that their method could achieve passenger costs with minimal differences within a reasonable timeframe, thereby validating the effectiveness of their approach.

**Table 5 Summary of passenger recovery studies using common recovery strategies**

| Studies | Disruptions | | | Recovery strategies | | Objective functions Min | | Solution methods |
|---|---|---|---|---|---|---|---|---|
| | AOG | AC | FD | FD | FC | DT | TRC | |
| Acuna-Agost et al. (2015) | Y | Y | Y | Y | Y | | Y | NPA |
| McCarty and Cohn (2018) | | | Y | Y | | Y | Y | BD |
| Wang et al. (2020) | Y | | | Y | Y | | Y | B&P+CG |

Notes: NPA: network pruning approach; BD: Benders decomposition; CG: column generation.

### 3.1.1.4 Aircraft and passenger recovery

The studies in the previous three sections have predominantly focused on individual recovery (aircraft, crew, passenger) and have developed relatively mature methods. However, as studies progress, the concept of integrated recovery emerged as a research hotspot. Furthermore, it is important to recognize that the optimal solution for individual recovery may not necessarily lead to the optimal solution for integrated recovery (Maher, 2016). In integrated recovery problems, the most prominent focus of scholars is the Aircraft and Passenger Recovery Problem (APRP). Finding a balance between the recovery costs for the airline and the losses for passengers is a debatable matter and has been a continual focal point in the research on APRP. This section only considered studies that explores APRP using common recovery strategies.

One of the notable instances on APRP is the challenge presented by the prominent European Operational Research Society (named ROADEF) in 2009. This challenge prompted numerous scholars to conduct studies and publish relevant papers based on this case (Bisaillon et al., 2010; Eggenberg et al., 2010; Jafari and Zegordi, 2011; Mansi et al., 2012; Jozefowiez et al., 2013; Sinclair et al., 2014, 2016; Zhang et al., 2016). These studies mainly focused on ARPR under disruptions of AOG, AC, and FD, while Zhang et al. (2016) consider airport flow control instead of AC, and Jafari and Zegordi (2011) ignored AC. In addition, flight delay and cancellation are most used recovery strategies in these studies. Also, total recovery cost is most commonly-used objective, while delay time is seldom used (Mansi et al., 2012). Finally, in terms of solution methods and tools, there is a wide variety, including Lingo (Jafari and Zegordi, 2011), CG (Eggenberg et al., 2010), and simple heuristic (Bisaillon et al., 2010; Mansi et al., 2012; Jozefowiez et al., 2013; Sinclair et al., 2014), and complex hybrid heuristic algorithms (Sinclair et al., 2016; Zhang et al., 2016).

In addition to the studies conducted on the ROADEF 2009 challenge, many other studies have been conducted on APRP (Jafari and Zegordi, 2010; Hu et al., 2015, 2016; Ng et al., 2020; Lee et al., 2021; Yan and Chen, 2021; Wandelt et al., 2023). Most of these studies are also based on the AOG disruption. Jafari and Zegordi (2010) demonstrated the superiority of their method compared to the baseline method that did not consider passengers. Hu et al. (2015) considered passenger reassignment and minimized recovery costs for both passengers and airlines, using a time-space network approach to solve the problem. In order to obtain high quality solutions within a short time, Hu et al. (2016) designed a GRASP algorithm to solve the model. Compared with the separate recovery method, their approach could achieve lower recovery costs in a shorter timeframe. Similar to the approach taken by Sinclair et al. (2014), Ng et al. (2020) developed an improved LNS algorithm to enhance running time and solution quality. There are also studies that consider other types of disruptions. For instance, Yan and Chen (2021) examined APRP in the context of a typhoon scenario. They employed a divide-and-conquer method in conjunction with CPLEX to solve the problem. Similarly, Lee et al. (2021) investigated the impact of climate change on flight disruptions and incorporated the airport capacity restrictions caused by temperature into their model. Building upon passenger recovery models developed by Bratu and Barnhart (2006) and Marla et al. (2017), as well as the aircraft recovery model from Lee et al. (2020), they enhanced and combined these models to construct a new APRP model and improved the dynamic RH algorithm to solve the problem. Results indicated that based on the RCP8.5 climate change scenario, the average daily total recovery costs for airlines increased by 15.7% to 49.4%. Wandelt et al. (2023) used an airline recovery model as a benchmark to evaluate the robustness of the air network. Their model minimized total recovery costs using a Hybrid Variable Neighborhood Search (HVNS) algorithm. Their findings indicated that for relatively simple disruption scenarios, complex network metrics could accurately quantify the impact. However, in more complex disruption scenarios, estimates derived from complex network metrics often underestimate the impact of node disruptions. This study presents a novel approach to examining both the airline recovery problem and the robustness of the air network, providing valuable insights for scholars in the field.

  Table 6 presents a summary of aircraft and passenger recovery studies using common recovery strategies. It can be seen that the studies on APRP have been comprehensive, with a significant focus on disruptions such as AOG and AC. Regarding recovery strategies, the primary emphasis remains on common strategies such as flight delay, cancellation, and aircraft swap. Only a few studies consider calling up reserve aircraft or ferrying aircraft strategies (Jafari and Zegordi, 2010, 2011). Regarding the objective functions, the primary aim has been to minimize the total recovery costs. The innovation in these studies primarily lies in the solution methods employed, ranging from the utilization of commercial solvers to the development of heuristic algorithms that generate faster and higher-quality solutions. However, it is essential to note that the mentioned studies predominantly focus on the recovery costs associated with passenger reassignments, delays, and cancellations, neglecting the heterogeneity in passenger decision-making and preferences when facing disruptions. This limitation could potentially deviate from real-world passenger behaviors. Addressing passenger decision-making during disruptions could be a future direction in the study of APRP.

Table 6 Summary of aircraft and passenger recovery studies using common recovery strategies

| Studies | Disruptions | | | | | | Recovery strategies | | | | | Objective functions | | | Solution methods/tools |
|---|---|---|---|---|---|---|---|---|---|---|---|---|---|---|---|
| | | | | | | | | | | | | Min | | Max | |
| | AOG | AC | FD | FC | AW | Others | FD | FC | SW | Ferry | Others | DT | TRC | Others | |
| Bisaillon et al. (2010) | Y | Y | Y | | | | Y | Y | | | | | Y | | LNS |
| Eggenberg et al. (2010) | Y | Y | Y | | | | Y | Y | Y | | | | Y | | CG |
| Jafari and Zegordi (2010) | Y | | Y | | | | Y | Y | Y | Y | Y | | Y | | LINGO |
| Jafari and Zegordi (2011) | Y | | Y | | | | Y | Y | Y | Y | | | Y | | LINGO |
| Mansi et al. (2012) | Y | Y | | | | | Y | Y | | | | Y | | Y | Oscillation-based algorithm |
| Jozefowiez et al. (2013) | Y | Y | Y | | | | Y | Y | | | | | Y | | NCF heuristic |
| Sinclair et al. (2014) | Y | Y | | | | | Y | Y | | | | | Y | | Improved LNS |
| Hu et al. (2015) | Y | | | | | | Y | Y | Y | | | | Y | | CPLEX |
| Hu et al. (2016) | Y | | Y | | | | Y | Y | Y | | | | Y | | GRASP |
| Zhang et al. (2016) | Y | | Y | | | Y | Y | Y | Y | | | | Y | | A three-stage math-heuristic |
| Sinclair et al. (2016) | Y | Y | | | | | Y | Y | | | | | Y | | LNS + CG |
| Ng et al. (2020) | Y | Y | | | | | Y | Y | | | | | Y | | Improved LNS |
| Lee et al. (2021) | | | | | | Y | Y | Y | Y | | | | Y | | RH |
| Yan and Chen (2021) | | | | | Y | | Y | Y | Y | Y | | | Y | | CPLEX |
| Wandelt et al. (2023) | | | | | | Y | Y | Y | | | Y | | Y | | HVNS |

Notes: LNS: large neighborhood search; NCF: new connections and flights; HVNS: hybrid variable neighborhood search.

**3.1.1.5 Aircraft and crew recovery**

Aircraft and Crew Recovery Problem (ACRP) is also one of the hot topics that scholars are paying attention to. Abdelghany et al. (2008) developed a recovery tool, DSTAR and aimed at minimizing resources assignment cost, total delay cost, and cancellation cost. They used a RH modeling framework with a greedy optimization strategy for solving. Maher (2016) formulated a model to minimize flight delay and cancellation costs, reserve crews costs, extra task costs, and crew deadhead costs. They designed an improved Column-and-Row Generation (CRG) algorithm to solve this problem. Results indicated that CRG obtained higher quality solutions compared to CG. Khiabani et al. (2022) considered disruptions in crew members with the objective of minimizing recovery costs, employing the BD algorithm to solve the problem. Their method demonstrated a reduction of 75.3% in delay time, a decrease of 21% in delayed flights, and a 75% reduction in cancelled flights compared to the string-based approach. Le and Wu (2013) studied ACRP under AOG and crew absence situations and designed an iterative tree growing with node combination method to solve this problem. Zhang et al. (2015) focused on AC and AOG scenarios, designing the TSH algorithm. Compared with a sequential algorithm and the RH algorithm proposed by Abdelghany et al. (2008), their algorithm generated high-quality solutions for all instances within 2 minutes, meeting practical operational requirements. Eshkevari et al. (2023) considered the disruption of flight delay and used a tabu search algorithm to solve the two-objective model. Similarly, Liu et al. (2023b) considered the impact of crew members' long connections on recovery costs and balanced the increase in swap or delay costs and the potential decrease in the total recovery cost by saving available working hours during long connections. Chen et al. (2020) explored the extreme rain situation in the ACRP and minimized the pair number objective, non-home base objective, non-short-connect objective, total delayed flights, and maximum delay time. An NSGA-II variant algorithm was used to solve the problem, showing improved performance compared to the baseline solution.

A summary of aircraft and crew recovery studies using common recovery strategies is presented in Table 7. It can be seen that the frequently employed strategies are flight delay and deadhead, followed by flight cancellation and aircraft swap. It can also be seen that the primary objective of ACRP is to minimize recovery costs. In addition, existing studies primarily focus on improving algorithms to reduce running time and providing high-quality recovery schemes within a reasonable timeframe, aiming to minimize losses for airlines.

**Table 7 Summary of aircraft and crew recovery studies using common recovery strategies**

| Studies | Disruptions | | | | | Recovery strategies | | | | | | Objective functions Min | | Solution methods |
|---|---|---|---|---|---|---|---|---|---|---|---|---|---|---|
| | AOG | AC | FD | CA | AW | FD | FC | SW | DH | URC | Others | TRC | Others | |
| Abdelghany et al. (2008) | | Y | | | | Y | Y | Y | Y | Y | | Y | | RH with greedy optimization strategy |
| Le and Wu (2013) | Y | | | Y | | Y | Y | | Y | | | Y | | An iterative tree growing with node combination method |
| Zhang et al. (2015) | Y | Y | | | | Y | Y | Y | Y | | | Y | | TSH |
| Maher (2016) | | Y | | | | Y | Y | | Y | | Y | Y | | CRG |
| Chen et al. (2020) | | | | Y | | Y | | Y | Y | Y | | | Y | NSGA-II variant |
| Khiabani et al. (2022) | | | | Y | | Y | Y | Y | Y | | | Y | | BD |
| Eshkevari et al. (2023) | | | Y | | | Y | Y | | Y | | Y | Y | Y | TS |
| Liu et al. (2023b) | | | Y | | | Y | Y | Y | Y | | | Y | | CG |

Notes: CRG: column-and-row generation; TS: tabu search.

**3.1.1.6 Integrated airline recovery**

The main participants in the airline recovery problem are threefold: aircraft, crew, and passenger. As the study framework continues to mature, some scholars have also started studying the Integrated Airline Recovery Problem (IARP). A list of related studies is presented Table 8. To the best of our knowledge, Lettovský (1997) first studied IARP in his doctoral dissertation and developed a model with the objective of maximizing airline profits. However, recognizing the complexity that arises by linking aircraft routes, crew assignments, and passenger flows, even minor disruptions could lead to highly challenging problem-solving situations. Hence, the problem was decomposed using a BD algorithm, defining the master problem as the aircraft schedule recovery problem, while the sub-problems were set as ARP, CRP, and PRP. This laid the theoretical groundwork for subsequent studies on IARP. Bratu and Barnhart (2006) also attempted to study IARP but did not account for crew recovery costs, strictly falling into APRP model. Petersen et al. (2012) comprehensively addressed IARP, applying theoretical methods to actual operations. They modified the model of Lettovský (1997) and improved the BD algorithm, considering constraint generation and column generation when dealing with crew recovery model. The proposed method was compared with the sequential approach in various disruption scenarios such as different levels of traffic flow control and AC. Their results significantly reduced passenger delays, flight disruptions, and various recovery costs, marking a milestone in the systematic study of IARP. Maher (2015) introduced cancellation variables to describe flight cancellations and passenger reassignment to alternative flights in the integrated recovery problem. They developed the model and designed CRG for solving. Results showed that the proposed method, on average, saved 22.98% in recovery costs. Evler et al. (2022) integrated aircraft turnaround into the whole recovery system, formulating the problem equivalent to the heterogeneous vehicle routing problem with time windows to minimize recovery costs. They used the RH dynamic algorithm for solving. Based on Frankfurt Airport operational data, their approach demonstrated a reduction in recovery costs compared to aircraft recovery, enhancing the airport resilience.

In the study of airline recovery problem, epidemics are among the disruption factors that must be considered. This is primarily due to the unprecedented negative impacts of diseases, such as SARS (Anand et al., 2003), H1N1(Bautista et al., 2009), Ebola virus (Bogoch et al., 2015), and most notably the globally pervasive COVID-19 (Huang et al., 2020), on the aviation industry (Li and Wang, 2023). Consequently, strategies need to be formulated during an epidemic to ensure that airlines can restore operations as swiftly as possible to sustain their existence. Xu et al. (2023a) are pioneers who have studied the airline recovery problem in the epidemic condition. They developed the Pandemic Airline Integrated Recovery (PAIR) model based on the disruption caused by COVID-19 among crews and passengers. The model aimed to minimize flight delay and cancellation costs, crew deadhead costs, costs incurred due to unallocated passengers, and personnel infections. The model was formulated as a mixed integer programming model and solved using a Branch-and-Cut (B&C) and LNS based solution framework. Based on case studies from the COVID-19 pandemic, the proposed method could obtain feasible solutions within reasonable time and outperformed deterministic methods and robust optimization methods in terms of recovery costs and personnel infections. The study outcome has opened up a new research direction for the management of airline disruptions during an epidemic.

Since using exact methods to solve the airline recovery problem is time-consuming and operationally infeasible in practice, studies are constantly exploring more efficient heuristic methods to address optimization problems. Machine Learning (ML) focuses on developing computational algorithms that enable computers to improve automatically through experience. Its powerful capabilities have made it a favored method among scholars (Jordan and Mitchell, 2015). Eikelenboom and Santos (2023) attempted to employ a ML ranker algorithm in addressing the IARP. They developed a model that aimed to minimize recovery costs and used actual operational data to train the model. They validated the effectiveness of the method using another set of data. Results indicated that employing ML reduced computational complexity and exhibited excellent performance in terms of running time and solution quality.

In summary, IARP is an exceedingly comprehensive problem that considers the recovery scenarios of all flight participants. However, its computational complexity is notably high, and traditional solution methods often fall short in providing satisfactory solutions. Consequently, there is a necessity to design more effective algorithms to tackle this challenge. Currently, there is a lack of research findings addressing this problem, and the development direction is quite limited. In the future, research efforts are expected to concentrate on enhancing and refining solution methods for IARP.

Table 8 Summary of integrated airline recovery studies using common recovery strategies

| Studies | Disruptions | | | | Recovery strategies | | | | | | Objective functions Min | | Solution methods |
|---|---|---|---|---|---|---|---|---|---|---|---|---|---|
| | AOG | AC | FD | Others | FD | FC | SW | DH | URC | Others | TRC | Others | |
| Petersen et al. (2012) | | Y | | Y | Y | Y | | | | | Y | | BD |
| Maher (2015) | | Y | | | Y | Y | | Y | Y | Y | Y | | CRG |
| Evler et al. (2022) | Y | | Y | | Y | Y | Y | | | | Y | | RH |
| Xu et al. (2023a) | | | | Y | Y | Y | Y | | | Y | Y | Y | B&C + LNS |
| Eikelenboom and Santos (2023) | | | | Y | Y | Y | Y | Y | Y | | Y | | ML |

Notes: B&C: branch-and-cut; ML: machine learning.

**3.1.1.7 Summary**

In the studies on common recovery strategies, the initial focus was on ARP, CRP, and PRP. However, studies on PRP are considered of little significance, and it is unlikely to be extensively studied in the future. In the realm of integrated recovery studies, the primary focus of scholarly attention lies on ACRP and APRP. Given that aircraft serve as carriers for both passengers and crews, any study that seeks to establish a connection between passengers and crews inherently involves aircraft participation. Consequently, there is a fundamental absence of dedicated integrated recovery studies specifically targeting passengers and crews, but rather a concentration on the aforementioned IARP. Commonly discussed disruptions in research contexts involve AOG and AC. Regarding recovery strategies, the most commonly-used strategies in studies encompass flight delays, flight cancellations, aircraft and crew swaps, crew deadhead, and the use of reserve crews. In terms of optimization objectives, there were initially multiple types, including minimizing delay time, recovery costs, and the number of delayed and cancelled flights. However, as studies have advanced, there has been a convergence in unifying all evaluation metrics into recovery costs. Consequently, minimizing the total recovery cost became the unified objective. In the realm of solution methods, initial models were relatively simple, so commercial solvers were predominantly used. As models become more complex, commercial solvers are unable to produce exact solutions. Hence, heuristic algorithms have been continuously improved to enhance solution quality and reduce running time. However, such studies merely represent an initial exploration of airline recovery problem. A more in-depth, comprehensive, and efficient exploration of this problem demands consideration of various factors and necessitates greater efforts in environmental conditions and research subjects.

**3.1.2 Aircraft cruise speed control strategy**

In actual operations, if a flight is already delayed due to a late takeoff, managers have the ability to instruct the respective pilot to increase the speed of the aircraft during the flight, provided that all necessary requirements are met. By doing so, they can reduce the flight time, thereby minimizing the spread of delays and enabling the flight to reach its destination as early as possible. This acceleration aims to minimize the impact of propagated delays. In many studies, the initial flight time, used as input, is not fixed. Instead, its variability is primarily achieved through the control of cruise speed. Aircraft Cruise Speed Control (ACSC) is widely used in the aviation industry for various purposes such as delay management, air traffic control, fleet route assignment and crew schedule (Cook et al., 2009; Bertsimas et al., 2011; Gürkan et al., 2016; Safak et al., 2017; Wen et al., 2020). Building on this, scholars have applied the ACSC strategy to address flight disruption. A summary of the relevant studies is shown in Table 9. Aktürk et al. (2014) initially incorporated cruise speed as a decision variable along with environmental constraints and cost coefficients into ARP. They proposed an enhanced quadratic programming formula to reduce computational complexity, aiming to manage disruption at lower costs. Lee et al. (2020) creatively introduced a combined reactive and proactive approach to airline disruption management. They optimized recovery decisions based on realized disruptions and predicted future disruptions. Using a Stochastic Reactive and Proactive Disruption Management (SRPDM) model and the RH algorithm, they minimized recovery costs by integrating the airport congestion queuing model. Results showed that the proposed method reduced recovery costs by 1-2% compared to static benchmark solutions. In terms of APRP, Arıkan et al. (2013) developed a Conic Quadratic Mixed Integer Programming (CQMIP) model with the objective of minimizing recovery costs and solved it using CPLEX. Results indicated that this method effectively handled disruptions occurring within airlines, providing an optimal trade-off between operational costs and passenger-related costs. Marla et al. (2017) developed an alternative approximate model based on a time-space network. This model neglected passenger delays and instead, focused solely on passenger disruption. Results showed that their method could potentially reduce passenger disruptions by up to 83% and save up to 5.9% of the airline's total costs. Yetimoğlu and Aktürk (2021) included cruise speed as a decision variable and formulated a model with the objective of maximizing airline profits. They proposed a TSH algorithm to solve the problem. Their integrated recovery strategy combining ACSC and aircraft swaps resulted in the lowest recovery cost. When it comes to IARP, Arıkan et al. (2017) conducted a CQMIP model that demonstrated the ability to solve instances in a reasonable amount of time. Building upon this study, Ding et al. (2023) expanded the work by developing the Proximal Policy Optimization-Variable Neighborhood Search (PPO-VNS) algorithm for solving IARP. Results indicated that the algorithm significantly reduced running time and demonstrated excellent solving performance.

In summary, the ACSC strategy can effectively absorb ground delays and propagated delays caused by preceding flights. However, it is important to note that while it can reduce delays and ensure passenger satisfaction, it may result in additional carbon emission costs for the airline. Therefore, finding a balance between ACSC and carbon emissions, as well as striking an equilibrium between airline costs and passenger-related costs, should be taken into consideration.

**Table 9 Summary of passenger-oriented airline recovery studies considering ACSC recovery strategy**

| Studies | Problem types | Objective functions | Solution methods/tools | Data size | | |
|---|---|---|---|---|---|---|
| | | | | # of Aircraft | # of Fleets | # of Flights |
| Aktürk et al. (2014) | ARP | Minimize TRC | CPLEX | 60 | 6 | 207 |
| Lee et al. (2020) | | | RH | N | 3 | 852 |
| Arıkan et al. (2013) | APRP | Minimize TRC | CPLEX | N | 6 | 1429 |
| Marla et al. (2017) | | | Xpress | N | N | 250 |
| Yetimoğlu and Aktürk (2021) | | Maximize profit | TSH | 53 | 6 | 208 |
| Arıkan et al. (2017) | IARP | Minimize TRC | CPLEX | 402 | N | 1254 |
| Ding et al. (2023) | | | PPO-VNS | 30 | 3 | 94 |

Notes: PPO-VNS: proximal policy optimization-variable neighborhood search; N: not included or mentioned.

### 3.1.3 Flexible aircraft maintenance strategy

The general idea behind the flexible aircraft maintenance strategy involves categorizing aircraft maintenance tasks into two types: fixed maintenance tasks and flexible maintenance tasks. For flexible maintenance tasks, strategies such as maintenance delay, maintenance cancellation, and maintenance swap can be employed. Among these, the aircraft maintenance swap strategy is considered to be the most innovative. Liang et al. (2018) first introduced this strategy, which focuses on minimizing the number of flight cancellations to avoid the associated high recovery costs. The general process is depicted in Fig.7. If Flight 2, operated by the green aircraft, experiences a delay without implementing the aircraft maintenance swap strategy, as illustrated in Fig.7(a), the scheduled delay cannot be compensated. Consequently, the flight would have to be cancelled, and the aircraft would remain at Airport B, until an appropriate maintenance time becomes available. However, if the maintenance swap strategy is executed when the delay occurs, as shown in Fig.7(b), the green aircraft can undergo maintenance at Airport A. This allows the flight operated by the green aircraft to continue without cancellation. Leveraging the flexibility of maintenance tasks can significantly reduce recovery costs.

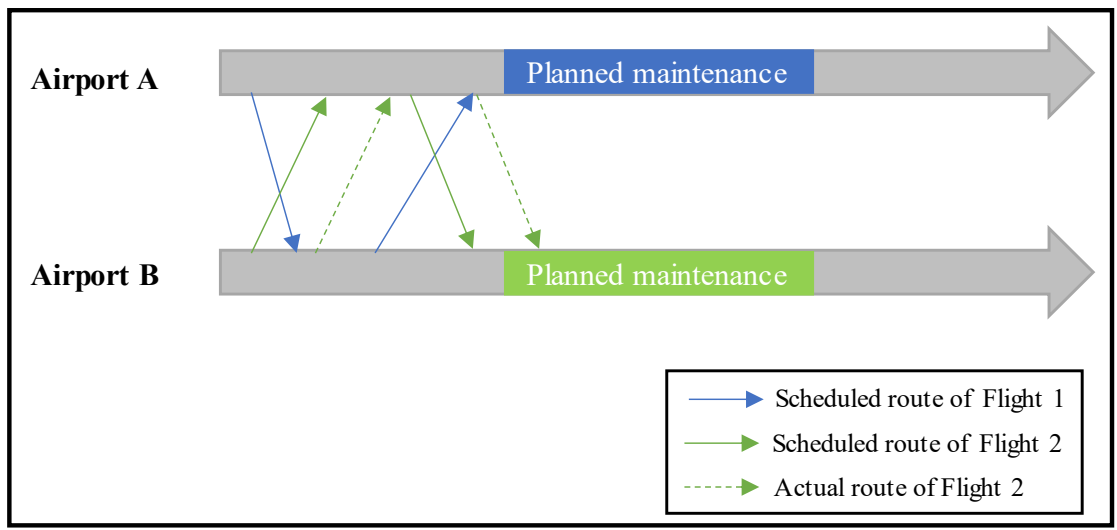

**(a) No maintenance swap**

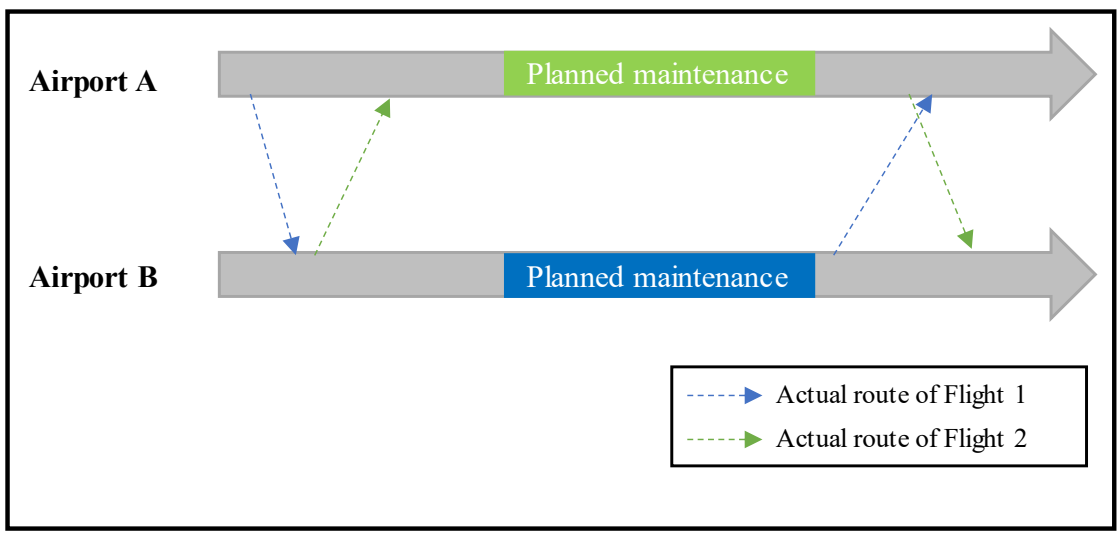

(b) Maintenance swap

Fig.7 Aircraft maintenance swap strategy

A summary of the existing studies on the flexible aircraft maintenance strategy can be found in Table 10. The table reveals that this particular strategy is predominantly adopted in the studies on ARP. As pioneers in this area, Liang et al. (2018) not only considered the flexibility of aircraft maintenance but also integrated airport capacity into the corresponding constraints. The objective of their study was to minimize recovery costs. To achieve this, they developed a model and used an improved CG algorithm to solve it. Results showed that the running time for corresponding large-scale problems averaged 6 minutes, and using aircraft maintenance swap strategy could potentially reduce recovery costs by up to 60%. The conclusions drawn by Liang et al. (2018) provided a novel perspective and study direction for subsequent studies considering the flexibility of aircraft maintenance. Li et al. (2022) improved the algorithm proposed by Liang et al. (2018), extending beyond the constraint that connections can only occur when the departure time of the first flight is earlier than that of the second flight to avoid missing the optimal solution. Additionally, after obtaining the optimal solution to the master problem, they employed the B&P framework to ensure optimal integer solutions. Huang et al. (2022) aimed to address the previous research gap regarding the evaluation of the quality of flight copies. They introduced the iterative cost-driven copy generation algorithm to assess flight copies, limiting their generation quantity, reducing problem scale, and minimizing running time. Rashedi et al. (2023) employed the supervised machine learning algorithm to solve ARP. They compared their algorithm with the Global Optimal Solution, Local Search 1 and Local Search 2 from Vink et al. (2020), CG at all nodes and CG at root nodes based on Liang et al. (2018). The comparison verified the algorithm's superiority and robustness. In addition, Vink et al. (2020) studied dynamic APRP and developed a model with the objective of minimizing recovery costs. Results indicated that their proposed method provided optimal solutions in an average time of 22 seconds for all instances, representing only 4% of the running time of dynamic global solutions. Regarding recovery costs, static global solutions were on average 19% lower than dynamic global solutions, suggesting that static ARP underestimates the required recovery costs.

Overall, using flexible aircraft maintenance for recovery is indeed a highly innovative approach, and there are many studies that have extensively explored this concept. However, it is worth mentioning that the existing studies primarily focused on ARP, with only Vink et al. (2020) delving into APRP. Future endeavors need to integrate such strategies into more integrated recovery problems to assess their effectiveness. Moreover, a practical consideration is that flexible aircraft maintenance involves written permission from the authorities. These frequent permissions might result in audits and become cumbersome each time the disruption happens.

**Table 10 Summary of passenger-oriented airline recovery studies considering flexible aircraft maintenance strategy**

| Studies | Problem types | Objective functions | Solution methods | Data size | | |
|---|---|---|---|---|---|---|
| | | | | # of | # of | # of |

|  |  |  | Aircraft | Fleets | Flights |
|---|---|---|---|---|---|
| Liang et al. (2018) |  | Improved CG | 44 | N | 638 |
| Li et al. (2022) |  | Improved CG | 38 | 5 | 172 |
| Huang et al. (2022) | ARP | ICCG | 162 | 9 | 789 |
| Rashedi et al. (2023) | Minimize TRC | SML | 220 | 2 | 1211 |
| Vink et al. (2020) | APRP | ASA | 100 | 2 | 600 |

Notes: ICCG: iterative cost-driven copy generation; SML: supervised machine learning; ASA: aircraft selection algorithm.

### 3.1.4 Multi-modal transportation strategy

The concept of multi-modal transportation strategy refers to the ability of managers or operators to select alternative transportation modes to send passengers to their destination when disruptions or delays occur. Referring to Xu et al. (2023b), An example of this concept is shown in Fig.8. As described in Fig.8, if Transportation Mode 2 from City 1 to City 2 experiences a disruption, relying solely on this mode would result in passengers waiting in City 2 until the service is restored. Due to the unpredictable duration of the disruption, it might cause substantial delays or even cancellations, leading to a large number of stranded passengers and potential economic losses. However, by considering alternative transportation modes within City 2 (such as Transportation Modes 1 and 3), managers can transport stranded passengers to their destination using alternative transportation modes. Taking Transportation Mode 1 as an example, passengers stranded at Facility 6 in City 2 could be transferred to Facility 2 within the same city using Transportation Mode 1 and from there, they can be further transported to their final destination.

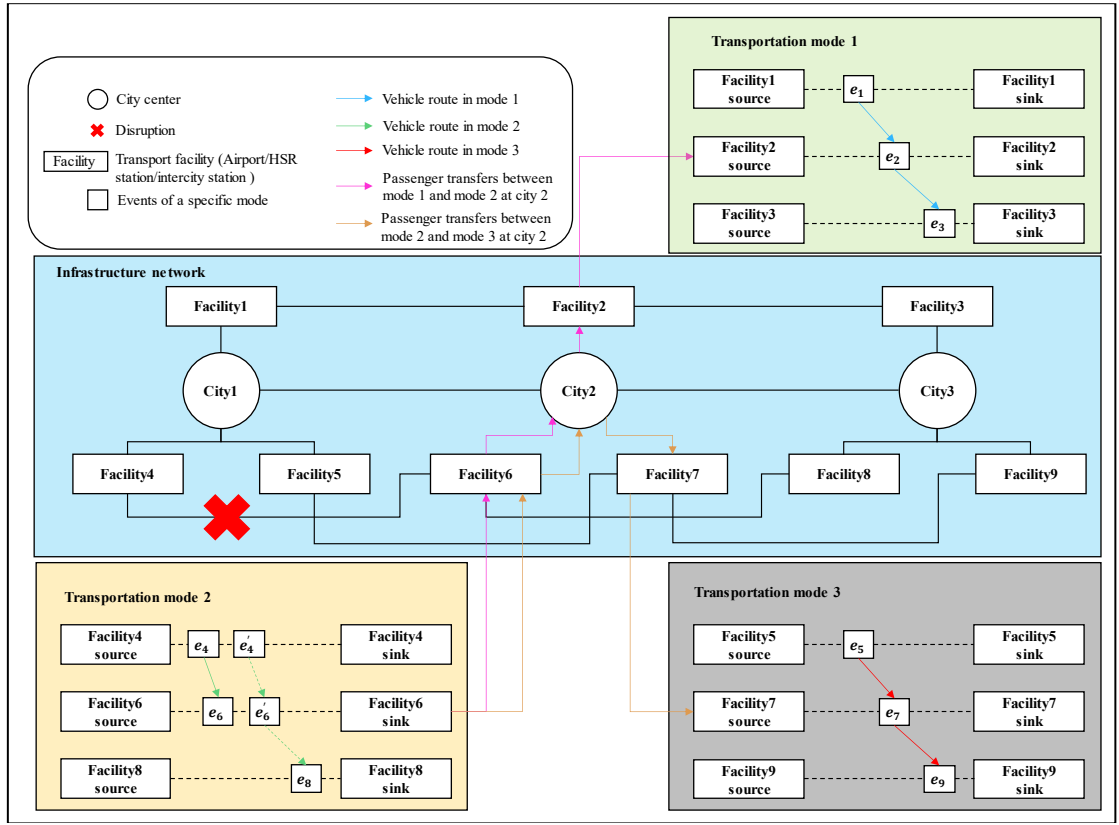

**Fig.8 Multi-modal transportation network**

Due to its flexibility, the multi-modal transportation strategy is frequently employed to recover from disruptions. One of the most common applications involves the design of substitute bus route under metro system disruption within cities (Kepaptsoglou and Karlaftis, 2010; Zhang et al., 2020). For intercity transportation, Liu et al. (2022a, 2023a) first conducted a study on the evacuation of stranded passengers during railway disruptions, utilizing a multi-modal transportation strategy. In the domain of airline disruption management, there has been a growing trend in the adoption of the multi-modal transportation strategy in recent years. Scholars have applied this strategy for various study purposes.

A summary of the existing studies on multi-modal transportation strategy can be found in Table 11. It can be seen that this strategy is predominantly explored in studies related to APRP. Zhang and Hansen (2008) were among the early scholars to consider employing such strategy to alleviate airport congestion, proposing the concept of Real-Time Intermodal Substitution (RTIMS). RTIMS uses ground transportation, primarily motor coaches, for airline recovery. Results revealed that the proposed method reduced recovery costs by 46% compared with the solution that do not consider RTIMS. Marzuoli et al. (2016), following an air crash, used ground transportation for airport recovery. They developed passenger-centric recovery schemes to balance costs and delays. Results showed that using motor coaches for passenger evacuation reduced the recovery costs for passengers, expediting airport recovery. Sun et al. (2022) introduced the concept of intermodalism to study APRP. They improved the Time-Band Network (TBN) model to reduce the running time, proving to be more efficient than other TBN-based models.

Moreover, several studies have been conducted on air-rail recovery strategy. Since High-Speed Rail (HSR) travel time is typically longer than air travel time, Liu et al. (2022b) introduced a resistance factor to illustrate the coordination costs between airlines and ground transport service providers. This factor reflected the airlines' willingness to financially contribute towards passenger transfers. A model was formulated with the objective of minimizing total recovery costs, validated using actual cases. Results indicated a 3.57% reduction in recovery costs for air-rail recovery compared to air-only recovery. Xu et al. (2023b) studied the Integrated Multimodal Mobility Under Network disRuptions (IMMUNER) problem, and developed an optimization model for air-rail recovery to minimize recovery costs. The model was verified using the Chinese multimodal transportation network. Results showed that flight cancellations had the most significant impact on passenger itineraries.

In summary, the existing studies on multi-modal transportation strategy for airline recovery have provided initial insights, establishing a solid foundation for further exploration of this problem. Most of existing studies in this field have predominantly focused on APRP. However, when it comes to passenger-related problems, these studies often assume that passengers will fully comply with the arrangements made by the airlines. Furthermore, most studies have not extensively considered passengers' subjectivity and heterogeneous factors (which will be extensively discussed in Section 3.1.5). Therefore, future study directions may include incorporating passengers' willingness and preferences into the study.

Table 11 Summary of passenger-oriented airline recovery studies considering multi-modal transportation strategy

| Studies | Problem types | Objective functions | Solution methods/ tools | Data size | | |
|---|---|---|---|---|---|---|
| | | | | # of Aircraft | # of Fleets | # of Flights |
| Zhang and Hansen (2008) | APRP | Minimize TRC | Approximation algorithm | N | 2 | 80 |
| Marzuoli et al. (2016) | | | Commercial solver | N | N | N |
| Sun et al. (2022) | | | MTBN, CPLEX | 188 | 13 | 628 |
| Xu et al. (2023b) | | | CRG | 198 | N | 454 |
| Liu et al. (2022b) | ARP | | PFCG, CPLEX | 118 | N | 295 |

Notes: MTBN: modified time-band network; PFCG: pseudo flight copies generation.

### 3.1.5 Passenger-centric recovery strategy

Most papers studying airline recovery problem overlook passengers' preferences and assume that passengers will fully comply with the arrangements made by airlines. However, in real situations, passengers facing flight disruptions tend to make different decisions. Existing papers highlight key decisions made by passengers, including waiting for the original flight (Yang and Hu, 2019; Hu et al., 2021b), endorsing to other flights (Yang and Hu, 2019; Hu et al., 2021b), transferring to HSR (Liu et al., 2022b;

Xu et al., 2023b), transferring to ground vehicles (Zhang and Hansen, 2008; Sun et al., 2022), and refunding tickets (Yang and Hu, 2019; Hu et al., 2021b). Passengers have various choices when faced with disruptions or changes to their travel plans. These choices can result in different outcomes including, successfully taking the flight to complete the itinerary as planned, being forced to refund due to willingness failure, opting to transfer to other transportation modes to complete the itinerary, and actively choosing to refund the ticket. The passenger-centric recovery strategy is illustrated in Fig.9.

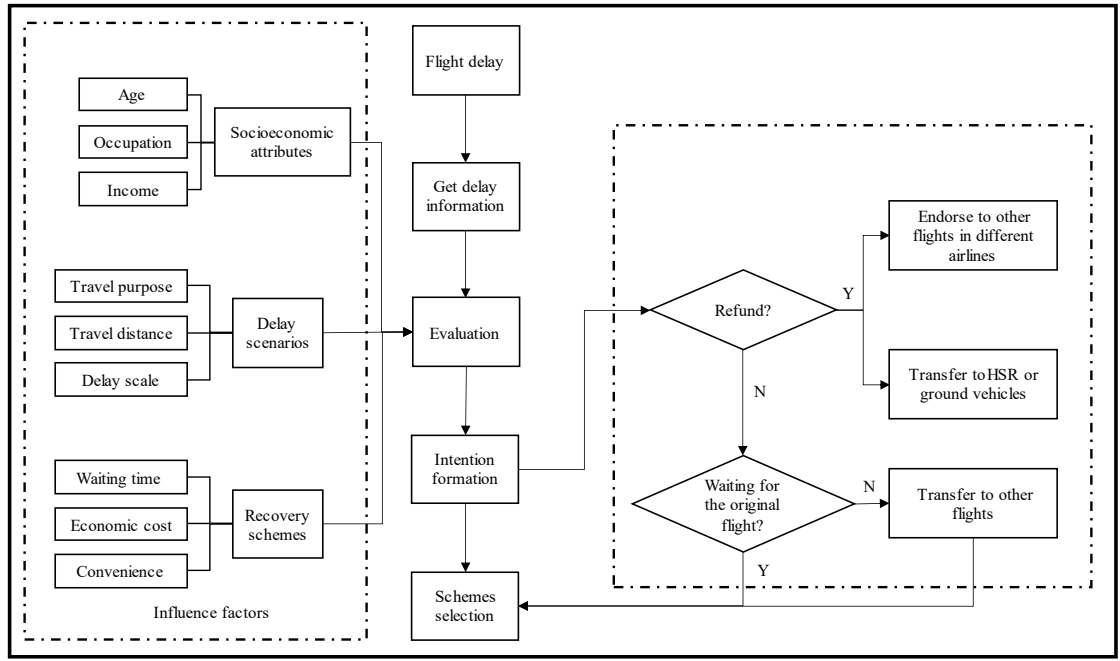

**Fig.9 Passenger-centric recovery strategy**

Passengers exhibit a strong subjectivity and possess distinct personal desires and preferences, which will greatly influence their decisions. Therefore, studies considering passenger preferences align more closely with real-life cases. Existing studies on this topic are listed in Table 12. Yang and Hu (2019) considered passengers' choices when facing flight delays, specifically the options to endorse to other flights or to refund their tickets. They formulated a bi-objective model aiming to minimize airline recovery costs and passenger utility loss, employing a loop-based multiobjective genetic algorithm for solving. Results showed that the approach decreased both recovery costs and passenger utility loss. Sensitivity analysis suggested that reducing turnaround time could decrease delays. Hu et al. (2021a) introduced prospect theory (Kahneman and Tversky, 1979) to refine passenger psychology and considered passengers' willingness to endorse and the willingness failure due to seat limitations, quantifying it as a loss in utility. In an effort to further study the problem, Hu et al. (2021b) provided a more detailed application of prospect theory in passenger choice behavior and designed a Multi-Directional and Stochastic Variable Neighborhood Search (MDSVNS) method, demonstrating lower costs and passenger utility loss compared to Multi-Directional Local Search heuristic (MDLS) and CPLEX. The conclusions provide a wealth of information and can serve as a valuable reference for future research. Cadarso and Vaze (2022) defined the "phantom passenger" concept for passengers who have confirmed reservations but fail to show up for their rebooked itinerary during disruptions. They formulated the IARP model, which minimizes total recovery costs by using CPLEX after linearization. Results indicated that their approach could restore airport operations at lower costs compared to neglecting phantom passengers. Chen et al. (2023) took into account the preferences of both the airline and passengers simultaneously and developed a model to minimize total airline recovery costs and passenger delay time and designed the Adaptive Non-dominated Sorting Genetic Algorithm-II based on Dominant Strengths (ANSGA2-DS) to solve the model. Results showed that their method outperformed the baseline methods, minimizing both airline delay costs and passenger delay time, with fare equity being optimal, reflecting the strength of their approach.

In summary, scholars have extensively studied the airline recovery problem, taking into account passenger preferences. Hu et al. (2021a, 2021b) have developed a relatively comprehensive theoretical framework, although primarily focusing on APRP. While Cadarso and Vaze (2022) extensively studied IARP, their consideration of psychological factors concerning passengers remained somewhat limited, neglecting scenarios where passengers might choose alternative transportation modes and the associated psychological impact. Consequently, future studies in this domain are expected to place greater emphasis on individual passenger psychological factors and conducting integrated recovery studies involving all stakeholders in flight operations.

Table 12 Summary of passenger-oriented airline recovery studies considering passenger-centric recovery strategy

| Studies | Problem types | Objective functions | Solution methods/ tools | Data size | | |
|---|---|---|---|---|---|---|
| | | | | # of Aircraft | # of Fleets | # of Flights |
| Yang and Hu (2019) | APRP | Minimize TRC and passenger utility loss | LBMGA | 59 | 1 | 209 |
| Hu et al. (2021a) | | | TBN and PTA | N | 10 | 414 |
| Hu et al. (2021b) | | | MDSVNS | 276 | 11 | 1038 |
| Chen et al. (2023) | | Minimize TRC and DT | ANSGA2-DS | N | N | 25 |
| Cadarso and Vaze (2022) | IARP | Minimize TRC | CPLEX | N | 5 | 1074 |

Notes: LBMGA: loop-based multiobjective genetic algorithm; TBN: time-band network; PTA: passenger transfer algorithm; MDSVNS: multi-directional and stochastic variable neighborhood search; ANSGA2-DS: adaptive non-dominated sorting genetic algorithm-II based on dominant strengths.

### 3.1.6 Clubbing of flights strategy

In addition to the aforementioned recovery strategies, scholars have also proposed the Clubbing of Flights (COF) strategy to address flight disruptions (Haider et al., 2023). The basic idea is described in Fig.10. Aircraft $a_1$ and $a_2$ operate flights $f_1$ and $f_2$, respectively, where Aircraft $a_2$ has a larger capacity than Aircraft $a_1$. Suppose Aircraft $a_1$ cannot operate Flight $f_1$ on time due to disruptions, without recovery strategies, Flight $f_1$ will be delayed. By implementing COF strategy, Flight $f_1$ and Flight $f_2$ are combined, and operated by Aircraft $a_2$, significantly reducing flight delays and recovery costs while ensuring passengers arrive at their destination as promptly as possible. Haider et al. (2023) innovatively adopted the COF strategy in conjunction with common recovery strategies to reduce flight delay costs. They acknowledged that maintenance swaps can introduce complex procedural audits, leading to the prohibition of such swaps. They developed a model with the objective of minimizing recovery costs and used the CG-LS algorithm to solve it. Results showed a significant 18.75% reduction in recovery costs by employing the COF strategy. Additionally, CG-LS significantly reduced running time, highlighting the efficiency and superiority of this method.

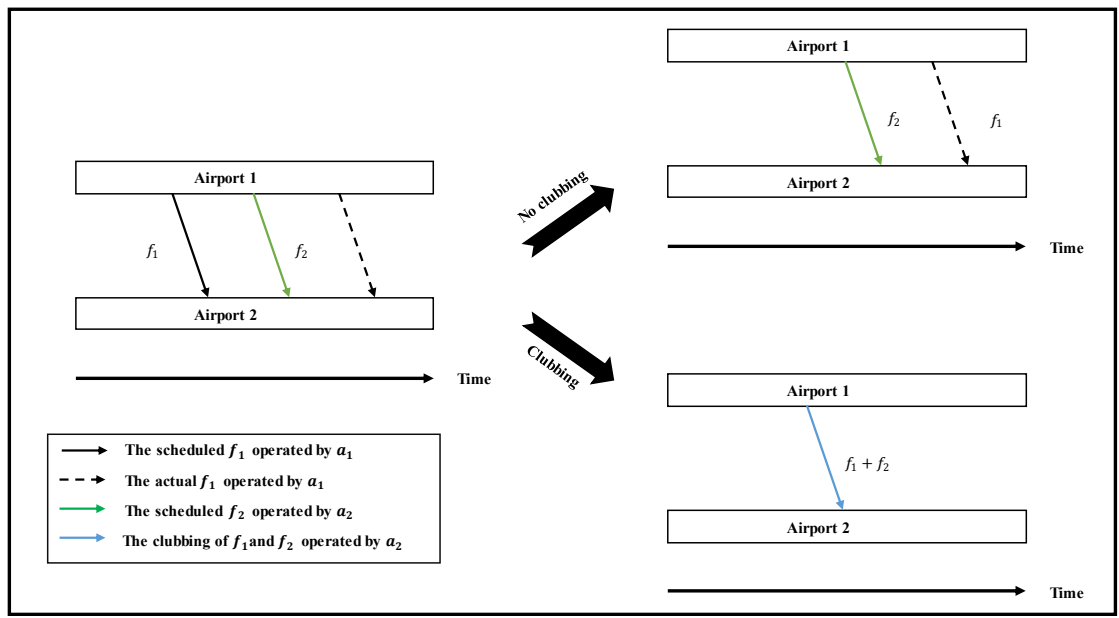

**Fig.10 The COF strategy**

The COF strategy offers a novel approach to address flight disruptions and has the potential to minimize recovery costs. However, its practical implementation requires consideration not only of the passenger numbers on both flights and the aircraft capacity, but also passengers' preferences. Additionally, comparing the effectiveness of this strategy with aircraft swap in airline recovery is another aspect that needs to be considered in future studies.

### 3.1.7 Summary

Studies on POARP have been extensive since 1984. Scholars have continually devoted themselves to exploring this problem, evolving from initial simple individual recovery problem to the current integrated recovery problem. The development has progressed from basic models to complex ones that represent various stakeholders' related costs or benefits. Initially, the constraints and assumptions were simplistic, but over time, they have become more detailed, reasonable, and realistic. This evolution is evident in the transition from simple case studies to complex examples that closely resemble real operational situations. Furthermore, the shift from simple solution methods to multiple high-efficiency solution algorithms that provide high-quality solutions demonstrates the continuous enrichment of the theoretical framework of POARP. The recommendations provided to managers are becoming increasingly valuable. In the following sections, efforts will be delved into disruptions, problem types, objective functions, and solution methods for a comprehensive summary.

### 3.1.7.1 Disruptions

Disruptions serve as the underlying causes for flights to be discontinued, prompting the evolution of studies on the airline recovery problem to address these challenges. Over the decades, numerous studies in the literature have explored airline recovery problem, focusing on different types of disruptions. As mentioned in Section 3.1.1.7, AOG and AC are the most widely studies disruption types. From the first paper addressing the airline recovery problem to present-day studies, a significant portion of the studies have been dedicated to investigating these two primary disruptions. Additionally, some studies have focused on examining the impact of adverse weather conditions on flights, aiming to offer reference plans for managers to effectively restore normal airport operations in the aftermath of weather-related disasters. Moreover, two disruptions that have gained attention in recent times are climate change and pandemic environments. Scholars have started exploring the implications of these disruptions on airline operations and recovery strategies. In the context of climate change, Lee et al. (2021) considered disruptions caused by global warming, which result in high temperatures that can impact smaller airports and hinder the takeoff of certain aircraft In regards to the pandemic environment, Xu et al. (2023a) addressed disruptions caused by COVID-19 infections among crews and passengers. These two disruptions are closely linked to current events, one stemming from the current global warming trend, and the other emerging from the recent COVID-19 pandemic that swept across the globe. These disruptions offer new pathways for future studies, enabling the extraction of real-life disruptions affecting flights from multiple perspectives. Scholars can then conduct more in-depth studies and develop relevant recovery plans to provide better advice to managers. This, in turn, enhances the resilience of the global aviation industry and enables a more composed and timely response to future disruptions.

### 3.1.7.2 Problem types

Regarding problem types, the distinction between individual and integrated recovery problems in the context of airline operations will be elaborated here. Individual recovery problems encompass three key areas: ARP, CRP, and PRP. Among these, studies on PRP are relatively limited. Since aircraft are the core of airport operations, ARP has attracted considerable attention and has accumulated a substantial theoretical framework concerning disruptions, recovery strategies, models, and solution methods. Considering that crews serve as both service personnel and operational staff, CRP has progressed from basic scheduling to a more comprehensive consideration of various factors (such as crew heterogeneity). Overall, individual recovery problem studies have become quite comprehensive, encompassing a wide range of considerations, with future study direction likely focusing on algorithm improvement.

Regarding integrated recovery problems, APRP stands out as a significant area of interest. This problem revolves around balancing aircraft recovery costs and passenger satisfaction, an area that scholars have long been committed to exploring. While there is a wealth of studies on this topic, a notable area that requires further development is to incorporate passenger subjectivity and decision-making into the model. In future studies, it would be beneficial to focus more on passenger heterogeneity and choice diversity. Passenger preferences should be considered in conjunction with flight-related costs to maintain the airline's reputation while ensuring passenger satisfaction on a foundation of minimal recovery costs. Similar to APRP, ACRP studies are extensive and rich in theory. Future studies in this field should concentrate on promoting fairness for crew members to maximize airline operational efficiency. By emphasizing the human aspect of crew scheduling recovery, it ensures that staff can serve passengers with a positive attitude, leading to improved customer satisfaction and overall operational performance. Given the triad of aircraft, crew, and passengers in airline recovery problem, the future study direction naturally aligns with IARP research, which has been well-established by several current publications. However, it is important to mention that IARP is large-scale, complex in modeling, challenging in solution, and requires efficient algorithms, posing a pressing need for solving these complexities. Optimal solutions for individual recovery problems might not be suitable for integrated recovery problems. Therefore, the challenge lies in identifying useful insights from individual recovery problems, improving and effectively incorporating them into the integrated recovery problem. In summary, with a clear distinction between individual and integrated recovery problems, future studies in the field of airline recovery should prioritize integrated recovery problems, staying updated with current developments, and providing timely theoretical support.

**3.1.7.3 Objective functions**

Basically, the airline recovery problem is an optimization problem, necessitating the establishment of mathematical programming models to achieve the desired or optimal values of relevant metrics, also known as objective functions. In the initial stage, studies primarily used the minimization of delay time as the objective to reduce flight delays and thereby expedite airport recovery. As studies progressed, strategies such as flight cancellations and aircraft swap became part of the considerations. Furthermore, as the scope of problem types expanded, it was essential to include metrics for crew members and passengers for further assessment. It then became evident that solely focusing on delay time was not sufficient to comprehensively represent the airport's recovery situation. As a result, it became necessary to select other objectives to capture the various dimensions of the recovery process. Some studies chose minimizing the number of disrupted passengers and the deviation of crew members from the original plans for evaluation as objectives. The ultimate objective of airline recovery problem is to reduce delays and ensure passenger satisfaction. To simplify the modeling process, numerous studies have converted crew and passenger-related metrics into costs and linked them with flight costs. The objective then becomes the minimization of TRC.

However, it is important to note that different studies have varied content within their objectives. While the general idea revolves around monetary costs, some metrics are challenging to quantify in terms of costs. Examples include passenger willingness or fairness in crew allocation. While some studies have considered these aspects in their objectives, they have not been widely adopted. Therefore, in future studies, incorporating certain human factors as metrics in the objective functions, in addition to minimizing costs, can better reflect universal conditions and make recovery schemes more convincing and inclusive.

**3.1.7.4 Solution methods**

Due to the increasing complexity and scale of the airline recovery problem, the solution methods have gone through several stages. Initially, simple models could be directly solved using commercial solvers, such as CPLEX. However, as the problem scales expanded, the efficiency of commercial solvers decreased, promoting scholars to explore other solution methods. As per the conclusions in Section 2.2, GA is one of the most commonly-used methods to address the airline recovery problem. GA, a renowned heuristic algorithm, was first introduced by Holland (1992). It avoids local optima through operations like crossover and mutation, enhances search efficiency, and is widely applicable in the domain of airline recovery. Moreover, enhanced versions of GA, such as the multi-objective genetic algorithm and NSGA-II (Deb et al., 2002), have been extensively applied in the field of airline recovery, yielding several representative research outcomes (Guo et al., 2005; Chen and Chou, 2017; Chen et al., 2023). Furthermore, due to airline recovery problem falling under NP-hard problem, CG has been widely used. This algorithm is an efficient method for large-scale linear optimization problems, originally proposed by Dantzig and Wolfe (1960). CG involves problem partitioning into master problem and subproblems, using iterations to solve, and exhibits high computational performance and good convergence. Studies developed CRG and its improved forms, which demonstrate efficiency and convergence for large-scale problems (Maher, 2015, 2016; Sinclair et al., 2016; Liang et al., 2018; Xu et al., 2023b). As the problem complexity increases, scholars have begun exploring other heuristic algorithms, ML, and various methods for solving. In summary, irrespective of the methods used, the aim is to obtain higher quality solutions within a reasonable time. Designing more efficient heuristic algorithms or combining heuristic algorithms with exact methods will be a future study direction.

## 3.2 Cargo-oriented airline recovery

Compared with POARP, there have been relatively few studies on the COARP. Existing studies primarily explores the service processes in air cargo operations among freight forwarders, airlines and terminal service providers (Feng et al., 2015), air cargo capacity allocation to maximize revenue (Delgado et al., 2019), and cargo planning problems, such as flight schedule (Derigs and Friederichs, 2012), and fleet route assignment (Xiao et al., 2022), rather than air cargo recovery. Due to the reduced flexibility and specific characteristics of air cargo transportation when it is placed in the belly of passenger aircraft, this study mainly focuses on COARP on freight aircraft.

A list of existing studies on COARP is presented in Table 13. To the best of our knowledge, Delgado et al. (2020) first studied Air Cargo Schedule Recovery Problem (ACSRP) and presented corresponding schemes based on three recovery policies: Scheme1 (Policy1: additional daily crews), Scheme2 (Policy2: additional crews between airports), and Scheme3 (Policy3: additional flights between airports per aircraft). The objective aimed to minimize recovery costs by establishing a Mixed Integer Linear Programming (MILP) model. They generated 24 instances based on four different disruptions to test the model by comparing it against the benchmark scheme. Results showed a cost reduction of 11.7%, 10.5%, and 9.8% for the three schemes compared with the benchmark scheme, confirming the effectiveness of the model. To further explore ACSRP and improve solution methods, Delgado and Mora (2021) set two strategies for transporting cargo: P1 (orders' cargo can be split for transport) and P2 (all items in an order must be transported together). They formulated a mathematical model for profit maximization and used an improved CG algorithm for solving. For method validation, they compared the results against two benchmark schemes. Results showed that both strategies outperformed the benchmark schemes, with P1 yielding the highest profit - increased by 54% and 15% compared with benchmark schemes. Although the profit of P2 was 3% lower than P1, its computational speed was significantly faster. Therefore, while separating cargo in an order could yield higher profits, it also complicates the problem, leading to longer running time.

Scholars have also studied integrated air cargo recovery problem. Mousavi et al. (2022) studied Aircraft and Air Cargo Recovery Problem (AACRP) and employed recovery strategies involving flight delay and cancellation, aircraft swap, and ACSC. They aimed to minimize recovery costs, for which they developed a MILP model and designed TSH algorithm for solving. When comparing the results with the benchmark plans, the operating costs were reduced by an average of 36%, along with an average reduction of 76% in running time. Huang et al. (2023) established two models for integrated cargo recovery problem: the arc-based model, named ACRP-A, and the string-based model, named ACRP-S. ACRP-A is relatively straightforward and can be directly solved using CPLEX. Meanwhile, ACRP-S is more complex, with authors using effective inequality constraints to transform the problem into ACRP-S*, and then designed a ML based CRG approach for solving. Results indicated that ACRP-S* could achieve lower recovery costs, with the ML-CRG algorithm saving 30% of the running time.

Table 13 Summary of cargo-oriented airline recovery studies

| Papers | Problem types | Objective functions | Solution methods/ tools | Main contributions |
|---|---|---|---|---|
| Delgado et al. (2020) | ACSRP | Minimize TRC | Gurobi | First study on ACSRP |
| Delgado and Mora (2021) | | Maximize profit | Improved CG | Develops improved CG that combines heuristics with mathematical programming techniques |

| Huang et al. (2023) | AACRP | Minimize TRC | ML-CRG | Uses ML technology |
| Mousavi et al. (2022) | | | TSH | Considers the fuel cost and develops a new algorithm |

Notes: ACSRP: air cargo schedule recovery problem; AACRP: aircraft and air cargo recovery problem; ML-CRG: machine learning based column-and-row generation.

## 3.3 Passenger-oriented versus cargo-oriented airline recovery

This section aims to provide a thorough summary of the connections and distinctions between studies of POARP and COARP. Both POARP and COARP are optimization problems requiring the establishment of mathematical models and the design of algorithms to solve them. Their objectives are aimed at reducing the recovery costs of airlines, ensuring customer satisfaction, and maintaining the airline's reputation. However, these two problems inherently possess certain differences, which are described from five aspects as presented in Table 14.

**(1) Timeliness**. In POARP, when passengers face disruptions, the primary purpose is to transport them to their destinations promptly, pacify their emotions, ensure their loyalty to the airline, and maintain the airline's reputation. However, in COARP, the main objective is to ensure that the cargo is delivered to the customer before the deadline without as much emphasis on the emotional aspects associated with customers.

**(2) Subjectivity**. POARP considers passenger autonomy in decision-making during disruptions, such as waiting, endorsing, using alternative transportation mode, or refunding, highlighting the need to consider passengers and aircraft at the same level. However, in COARP, cargo, being devoid of subjective will, follows ground staff schedules and airline arrangements.

**(3) Flexibility**. For POARP, efforts should be made to minimize passenger transfers to ensure a direct route to the destination, as multiple transfers can cause dissatisfaction. Conversely, in COARP, not having subjective preferences, cargo allows for more flexible routes, multiple transshipments, or segmented deliveries within the allowed timeframe.

**(4) Transferability**. In POARP, passengers have the independency to undertake transit operations, needing minimal ground facilities aside from shuttle buses and boarding bridges. However, in COARP, cargo relies on airline staff and numerous facilities during transshipment, stressing the allocation of airport resources.

**(5) Combinability**. Cargo can sometimes be transported in the belly of passenger aircraft, optimizing airline resources. Cargo flights, however, follow more flexible routes with multiple landings, not accommodating passenger arrangements onto cargo aircraft. Thus, part of the cargo can be combined in passenger flights for recovery purposes, but the reverse is not applicable.

In conclusion, while POARP and COARP share similarities, they fundamentally differ in various aspects. Therefore, in-depth analysis catering to their unique characteristics is essential when devising recovery schemes. This is crucial to ensure the efficiency of both passenger and cargo airlines, especially for airlines that handle both passenger and cargo, and to explore the trade-offs and timely recovery strategies between passenger and cargo, demanding further exploration by scholars in this field.

Table 14 Comparison between POARP and COARP

| Types | Comparison contents | | | | |
|---|---|---|---|---|---|
| | Timeliness | Subjectivity | Flexibility | Transferability | Combinability |
| POARP | As quickly as possible | √ | × | √ | × |
| COARP | Before deadline | × | √ | × | √ |

# 4 Conclusions and future study directions

This paper summarized and discussed the recent research progress on airline recovery problem. Firstly, relevant studies were searched using keywords, and a preliminary analysis was conducted using visualization tools to examine publications in terms of the number of papers, keywords, and sources. The discussion was then divided into POARP and COARP, with relatively fewer papers on COARP due to its recent attention. Regarding the studies on POARP, this paper took a different approach compared to existing review papers. It categorized the recovery strategies into common recovery strategies, ACSC strategy, flexible aircraft maintenance strategy, multi-modal transportation strategy, passenger-centric recovery strategy, and COF strategy. This paper provided a comprehensive summary of POARP, covering disruptions, problem types, objective functions, and solution methods. Finally, the paper discussed the differences between POARP and COARP in terms of timeliness, subjectivity, flexibility, transferability, and combinability.

After conducting a systematic review of airline recovery problem, the following recommendations for future studies are given.

**(1) Conduct studies on dynamic or real-time airline recovery.** When flights face disruptions, there is a clear need for swift recovery within a short timeframe. However, existing studies primarily focus on minimizing recovery costs and explore the recovery of predefined disruptions. It is noted that only a very few papers consider the importance of dynamic and real-time recovery (Lee et al., 2020; Vink et al., 2020), resulting in a relatively underdeveloped research system. Therefore, it is recommended to conduct studies on dynamic and real-time airline recovery. Inspiration from other transportation modes, such as the comparatively mature field of dynamic railway recovery (Zhan et al., 2015; Zhu and Goverde, 2021), can provide valuable insights for future airline recovery study.

**(2) Incorporate human factors into the modeling.** Regardless of POARP or COARP, a primary goal is to provide satisfactory service to customers to enhance their loyalty to airlines. However, current studies tend to prioritize airline recovery over passengers' individual factors. Addressing human factors is crucial. For POARP, an emphasis on passenger decision-making during disruptions and an analysis of their psychological states should be incorporated to quantify subjective factors for a more realistic approach. Moreover, in COARP, more consideration to crew factors in duty distribution, particularly regarding fairness, is essential. The air cargo and crew recovery problem can be a potential future study direction for cargo recovery problems.

**(3) Focus more on the application of multi-modal transportation strategy in airline recovery problem.** Implementing various transportation modes for airline recovery is a promising concept. However, existing studies using multiple transportation modes for airline recovery, such as RTIMS (Zhang and Hansen, 2008), are limited. Often, other transportation modes are assumed as virtual flight legs, neglecting their unique characteristics. Consideration of the itinerary's length and its relation to transportation mode choice is critical. For longer itineraries, using HSR for passenger transportation might be uneconomical. Additionally, negotiations between airlines and other operators regarding costs are essential aspects requiring further exploration. In future studies, a comprehensive consideration of the characteristics of multiple transportation modes is essential, aiming to better integrate them with airport operations. In terms of cargo transportation, maritime transportation serves as a primary alternative mode, while for passenger transportation, railways are a primary alternative mode. It is also crucial to incorporate passengers' willingness to choose alternative transportation modes into the scope of study on POARP.

**(4) Explore how other airport processes aid in the recovery of disrupted flights.** Airport processes are interconnected, and improving the efficiency of other processes may effectively counter flight disruptions. Evler et al. (2022) demonstrated that considering aircraft turnaround can significantly reduce recovery costs. Future studies could consider other processes, such as enhancing check-in efficiency or organizing orderly passengers queue at boarding gates, to explore their positive effects on airline recovery.

**(5) Combine the robust proactive scheduling and reactive airline recovery strategies to enhance airport resilience.** Airline disruption management can be categorized as proactive and reactive processes (Ogunsina et al., 2022). It is essential not only to study airline recovery (reactive) or robust scheduling (proactive) separately but also to integrate both approaches. This integration facilitates proactive (tactical) disruption management before execution, reactive (operational) disruption management during execution, and proactive (strategic) disruption management after execution, collectively enhancing the airport's resilience. Future studies should focus on how to extract relevant information from vast data, considering various stakeholders' interests, and integrating proactive and reactive strategies to better cope with uncertain disruptions.

**(6) Develop more efficient algorithms to optimization models for airline recovery problems.** One of the future study directions for all optimization problems is the improvement and validation of their solution algorithms. Airline recovery problem is no exception. While various methods have been employed to tackle this problem, it is important to acknowledge that there are still limitations in the algorithms used to solve it. For instance, exact algorithms often underperform in large-scale instances, and traditional heuristic algorithms, like GA and SA, have shortcomings. Some studies in the literature have used hybrid heuristic algorithms, yet their universal applicability requires further verification. Regarding ML, when using RL, the careful consideration of reward function setting to cater to diverse scenarios is crucial. Future studies should aim to design more efficient methods by amalgamating exact methods, heuristic algorithms, and ML to adapt to a broader range of disruption scenarios, rendering the solutions more universally applicable.